\documentclass{article}
\usepackage{amsmath}
\usepackage{amsfonts}
\usepackage{amssymb}
\usepackage{amsthm}

\def\a{\alpha}
\def\b{\beta}
\def\d{\delta}
\def\e{\varepsilon}
\def\g{\gamma}

\def\D{\Delta}

\def\O{\Omega}

\def\ot{\otimes}
\def\dag{\dagger} 
\def\no{\nonumber}

\def\S{\tilde{S}}

\def\X{\bar{X}}
\def\Y{\bar{Y}}
\def\Z{\bar{Z}} 

\def\R{{\cal R}}

\newcommand{\reff}[1]{(\ref{#1})}

\newtheorem{Theorem}{Theorem}
\newtheorem{Definition}{Definition}
\newtheorem*{Definition*}{Definition}
\newtheorem{Proposition}{Proposition}
\newtheorem{Lemma}{Lemma}
\newtheorem{Corollary}{Corollary}
\newtheorem*{Corollary*}{Corollary}

\theoremstyle{remark}
\newtheorem*{remark*}{Remark}

\begin{document}

\begin{titlepage}
\vskip.3in

\begin{center}
{\Large \bf  Quasi-Hopf $*$-Algebras}
\vskip.3in
{\large M.D. Gould and T. Lekatsas } 
\vskip.2in
{\em Department of Mathematics\\ The University of Queensland\\ Brisbane,
     Qld 4072, Australia

Email: tel@maths.uq.edu.au
\vskip.2in
\today}
\end{center}

\vskip 2cm

\begin{abstract}
We introduce quasi-Hopf $*$-algebras i.e. quasi-Hopf algebras equipped with a conjugation (star) operation. 
The definition of quasi-Hopf $*$-algebras proposed ensures that the class of quasi-Hopf $*$-algebras is
closed under twisting and additionally, that any Hopf $*$-algebra becomes a quasi-Hopf $*$-algebra via twisting. 
The basic properties of these algebras are developed. The relationship between the antipode and star structure
is investigated. Quasi-triangular quasi-Hopf $*$-algebras are introduced and studied.
\end{abstract}

\end{titlepage}
\vskip 3cm


\section{Introduction\label{intro}}

Many of the ideas and constructions from the theory of Hopf algebras have analogues in the quasi-Hopf algebra setting.
Examples include the quantum double construction~\cite{bc1,hn2,majid97,pan1}, the Tannaka-Krein theorem~\cite{majid92}, the existence of integrals~\cite{bc2,hn1,po}, the construction of link invariants~\cite{ac,Lin00} and extension to the superalgebra case~\cite{oqhsa,Zha98}, amongst others.
 
In the Hopf algebra setting, Hopf algebras that admit a conjugation or star operation are well known~\cite{Maj00, Sch93, Wor87} and are called $*$-Hopf algebras or Hopf $*$-algebras. The introduction of a $*$-structure is not part of the algebraic formalism of Hopf algebras but becomes necessary for physical applications, such as in quantum mechanics where unitarity is a requirement. The notion of a $*$-structure has been extended to the weak Hopf algebra case~\cite{bns} and also to the braided setting~\cite{majid94}, but the quasi-Hopf algebra case appears to have been neglected in the literature. Quasi-Hopf algebras have applications in conformal field theory~\cite{DPR, DVVV} and in the theory of integrable models (via elliptic quantum groups~\cite{Enr97,Fod94,Fel95,Fro97,Jim97,Zha98}). It is expected that $*$-structures should arise naturally
in such physical applications.

In this paper we introduce quasi-Hopf $*$-algebras ($*$-quasi-Hopf algebras, $*$-QHA). Our definition is motivated by the twisting construction
of Drinfeld~\cite{Dri90} which turns a Hopf algebra $H$ into a quasi-Hopf algebra.
The twisting operation changes the co-algebra structure of $H$ in such a way that the twisted co-product is no longer co-associative. The algebra
structure of $H$ is not affected by twisting.
The axioms introduced by Drinfeld for quasi-Hopf algebras ensure that any Hopf algebra will be twisted into a quasi-Hopf algebra and
that the class of quasi-Hopf algebras is itself closed under the twisting operation.
This larger class of algebras contains the Hopf algebras within it since every Hopf algebra is trivially a quasi-Hopf
algebra.  

A Hopf algebra $H$ may be equipped with a $*$-operation $\dagger:H \rightarrow H$, whenever the base field $\mathbb{F}$ over which it is defined admits a conjugation operation. A Hopf $*$-algebra~\cite{Wor87} is a Hopf algebra equipped with a $*$-operation $\dagger: H \rightarrow H$ such that on the algebra part of $H$, the conjugation $\dagger$ obeys the usual axioms of a $*$-algebra, and such that on the co-algebra part, the co-product $\D:H \rightarrow H \otimes H$ and the co-unit $\e: H \rightarrow {\mathbb F}$ are $*$-algebra homomorphisms.
The antipode $S$ of a Hopf $*$-algebra necessarily obeys $S(a)^{\dagger}=S^{-1}(a^{\dagger}), \forall a \in H$. This is a direct consequence of the uniqueness of the antipode.

Since a Hopf $*$-algebra is a Hopf algebra, twisting changes the Hopf algebra part into a quasi-Hopf algebra. As twisting does not affect the algebra structure of $H$, the $*$-algebra part of $H$ is unchanged. The $*$-structure on the co-algebra $H$ is twisted in such a way that the twisted co-product $\D_F$ is no longer a $*$-algebra homomorphism.
Nonetheless, $\D_F$ is a $*$-algebra homomorphism up to conjugation by the self adjoint twist 
$\O = (F F^{\dagger})^{-1}$
\begin{equation*}
\D_F(a)^{\dagger} = \O \D_F(a^\dagger)\O^{-1}, \quad \forall a \in H.
\end{equation*}
Twisting makes $H$ into a quasi-Hopf algebra and thus it has a co-associator $\Phi_F$ (induced by $F$). We show that the
co-associator $\Phi_F$ is related to its conjugate inverse $(\Phi^{\dagger}_F)^{-1}$ by the same twist $\Omega$ i.e.
\begin{equation*}
(\Phi^{\dagger}_F)^{-1} = (\O \otimes 1) \ (\D \otimes 1) \O \ \Phi_F \ (1 \otimes \D) \O^{-1} \ (1 \ot \O^{-1}).
\end{equation*}

Our definition of quasi-Hopf $*$-algebras is motivated by these observations. We define a $*$-quasi-Hopf algebra to be a quasi-Hopf algebra equipped with a conjugation $\dagger:H \rightarrow H$ and a twist $\O \in H \otimes H$ such that
\begin{equation*}
\e(a^{\dagger}) = \overline{\e(a)}, \quad \forall a \in H
\end{equation*}
\begin{equation*}
\D(a)^{\dagger}=\O \D(a^{\dagger})\O^{-1}, \quad \forall a \in H
\end{equation*} 
\begin{equation*}
(\Phi^{\dagger})^{-1} = (\O \otimes 1) \ (\D \otimes 1) \O \ \Phi \ (1 \otimes \D) \O^{-1} \ (1 \ot \O^{-1}) \equiv \Phi_{\O}.
\end{equation*}
This definition ensures that any Hopf $*$-algebra $H$ is twisted into a quasi-Hopf $*$-algebra. We show that the class of
quasi-Hopf $*$-algebras is closed under twisting. Unlike the Hopf algebra case, the antipode $S$ of a quasi-Hopf algebra is not unique. In the quasi-Hopf $*$-algebra setting this means that $S$ is not forced to satisfy any particular condition. 

We develop the general theory of quasi-Hopf $*$-algebras and investigate the relationship between the antipode $S$ and the conjugation operation $\dagger$ on $H$.  The effect of the Drinfeld twist on the $*$-canonical element $\O$ is determined and an explicit expression for the conjugate of the Drinfeld twist is derived. Quasi-triangular quasi-Hopf $*$-algebras are introduced. As in the Hopf algebra case, there are
two natural classes of quasi-triangular quasi-Hopf $*$-algebras. In the type I case, the $R$-matrix satisfies $(\R^{\dagger})^{-1}=\O^T \R \O^{-1}$, whilst for the type II case it satisfies $(\R^{\dagger})^{-1}=\O^T (\R^T)^{-1} \O^{-1}$. These reduce to the antireal and real cases of Majid~\cite{majid}, respectively in the Hopf algebra case where $\O = 1 \otimes 1$.
 
A further motivation for our definition comes from the quantised universal enveloping algebra $U_q(L)$ of a semi-simple Lie algebra $L$,
when $q \in \mathbb{C}$ is a complex phase. For $q$ real and positive $U_q(L)$ is a Hopf $*$-algebra. However, when $q \in \mathbb{C}$ is a complex phase, $\bar q = q^{-1}$ on conjugation, so that the conjugate of the co-product has the natural structure
of the opposite co-algebra i.e.
\begin{equation*}
\Delta(a)^{\dagger} = \Delta^T(a^\dagger), \forall a \in U_q(L).
\end{equation*}
Thus when $q$ is a phase, $H$ is not a Hopf $*$-algebra as noted in~\cite{Maj00}. Since $U_q(L)$ is
quasitriangular, it has an $R$-matrix $\R$. Now $\R$ is a twist and satisfies $\R \Delta(a) = \Delta^T (a) \R$ so that
\begin{equation*}
\Delta(a)^\dagger = \R \Delta(a^{\dagger}) \R^{-1}.
\end{equation*}
We take $U_q(L)$ to be a quasi-Hopf algebra with trivial co-associator $\Phi = 1 $. Now, $\Phi_{\R} = 1 \ot 1 \ot 1$ follows from the
quantum Yang-Baxter equation, so that the $(\Phi^{\dagger})^{-1} = \Phi_{\R}$ is trivially satisfied.
Thus $U_q(L)$ for $q$ a phase has the structure of a $*$-quasi-Hopf algebra with $*$-canonical element $\R$.


\section{Preliminaries}\label{prelim}

We begin by recalling the definitions and basic properties of quasi-bialgebras (QBA) and quasi-Hopf algebras (QHA).
\begin{Definition}
	A quasi-bialgebra $H$ is a unital associative algebra over a field $\mathbb{F}$,
	equipped with algebra homomorphisms $\e: H \rightarrow \mathbb{F}$ (co-unit),
	$\D: H \rightarrow H \ot H$ (co-product) and an invertible element 
	$\Phi \in H \ot H \ot H$ (co-associator), satisfying
	\begin{eqnarray}
		        (\e \ot 1) \D = & 1 & = (1 \ot \e) \D \label{counit}\\
		        (1 \ot \D) \D(a) & = & \Phi^{-1}(\D \ot 1) \D(a) \Phi,\quad \forall a \in H \label{q-co}     \\
                        (\D \ot 1 \ot 1) \Phi \ (1 \ot 1 \ot \D) \Phi & = &
		        (\Phi \ot 1) \ (1 \ot \D \ot 1) \Phi \ (1 \ot \Phi) \label{pentagon} \\  
	                (1 \ot \e \ot 1) \Phi & = & 1.                                                   \label{epsphi}   
	\end{eqnarray}
	A quasi-bialgebra $H$ equipped with an algebra anti-homomorphism $S: H \rightarrow H$ (antipode)
	and canonical elements $\a, \b \in H$ satisfying
	\begin{eqnarray}
	         \sum_\nu S(X_\nu) \a Y_\nu \b S(Z_\nu) & = &  1  = \sum_\nu \X_\nu \b S(\Y_\nu) \a \Z_\nu       \label{Sphi}\\
		 \sum_{(a)} S(a_{(1)}) \a a_{(2)} & = & \e(a) \a, \quad   \sum_{(a)} a_{(1)}\b S(a_{(2)}) = \e(a) \b, 
	                                                                                      \quad \forall a \in H \label{Sab}. 
	\end{eqnarray}
	is called a quasi-Hopf algebra. 
\end{Definition}

Above we have used Sweedler's~\cite{Swe} notation for the co-product 
\begin{eqnarray*}
	\D(a)=\sum_{(a)} a_{(1)}\ot a_{(2)}, \quad \forall a \in H.
\end{eqnarray*}
The co-product is no longer co-associative for QHA necessitating an extension to
Sweedler's notation 
\begin{eqnarray*}
	(1\ot\D)\D(a)  & = & \sum_{(a)} a_{(1)} \ot \D(a_{(2)})= \sum_{(a)} a_{(1)} \ot a_{(2)}^{(1)} \ot a_{(2)}^{(2)} \no \\
	(\D\ot 1)\D(a) & = & \sum_{(a)} \D(a_{(1)}) \ot a_{(2)}= \sum_{(a)} a_{(1)}^{(1)} \ot a_{(1)}^{(2)} \ot a_{(2)}. 
\end{eqnarray*}
For the co-associator we follow the notation of \cite{mdgtl,cas,oqhsa} and write
\begin{eqnarray*}
	\Phi      = \sum_\nu X_\nu \ot Y_\nu \ot Z_\nu~,\qquad 
	\Phi^{-1} = \sum_\nu \X_\nu \ot \Y_\nu \ot \Z_\nu.
\end{eqnarray*}
We adopt the above notation throughout and in general omit the summation sign from expressions,
with the convention that repeated indices are to be summed over.

It follows from equations \reff{counit}, \reff{pentagon} and \reff{epsphi} that the
co-associator $\Phi$ has the following useful properties
\begin{align*}
	(\e \ot 1 \ot 1)\Phi = 1 = (1 \ot 1 \ot \e) \Phi.
\end{align*}
Throughout we assume bijectivity of the antipode $S$ so that $S^{-1}$ exists.
The antipode equations~\reff{Sphi},~\reff{Sab} imply 
\begin{eqnarray*}
	\e(\a) \e(\b) & = & 1 \\
	     \e(S(a)) & = & \e(S^{-1}(a)) = \e(a),\quad \forall a \in H.
\end{eqnarray*}

Let $H$ be a QHA, an element $F \in H \ot H$ is called a twist (or gauge transformation) if it is invertible and satisfies
the co-unit property
\begin{equation}
	(\e \ot 1)F = (1 \ot \e)F = 1. \label{twist-counit}
\end{equation}
The operation of twisting operation allows one to construct a new QHA $H_F$ from $H$, called the twisted structure
induced by $F$, with the same antipode and co-unit, but with co-product, co-associator and canonical elements given by
\begin{eqnarray}
	\D_F(a) & = & F \D(a) F^{-1}, \quad \forall a \in H \label{TWDEL} \\
	\Phi_F  & = & (F \ot 1)(\D \ot 1)F \ \Phi \ (1 \ot \D)F^{-1} (1 \ot F^{-1}) \label{TWPHI}\\
	   \a_F & = & m \cdot (1 \ot \a)(S \ot 1)F^{-1}, \quad \b_F = m \cdot (1 \ot \b)(1 \ot S)F.       \label{TWCAN}
\end{eqnarray}
Above $m:H \ot H \rightarrow H$ is the multiplication map $m \cdot (a \ot b) = ab$.

Let $T:H \ot H \rightarrow H \ot H$ be the usual twist map $T(a \otimes b) = b \ot a$.
Recall that a quasi-Hopf algebra $H$ is also a quasi-Hopf algebra with the opposite co-product $T \cdot \D$ as follows,
\begin{Proposition}\label{P1}
	Let $H$ be a QHA. Then the opposite QHA, $H^{cop}$ is a QHA with co-product
	$\D^T = T\cdot \D $, co-associator $\Phi^T \equiv \Phi_{321}^{-1}$, antipode $S^{-1}$
	and canonical elements $\a^T \equiv S^{-1}(\a),~\b^T \equiv S^{-1}(\b)$.
\end{Proposition}


\section{Twisting on Hopf $*$-algebras}
In this and the following sections we take the base field to be the field of complex numbers $\mathbb{C}$.
Recall that a bi-algebra is a QBA with trivial co-associator $\Phi =1 \ot 1 \ot 1$. Similarly a Hopf algebra
is a QHA with trivial co-associator and trivial canonical elements $\a=\b =1$.

\begin{Definition}
	A bi-algebra $H$ is called a $*$-bi-algebra if it admits an antilinear map
	$\dagger: H \rightarrow H$ (conjugation operation) satisfying
	\begin{eqnarray} 
		(a^\dag)^\dag & = & a \label{ADAG} \\
		    (ab)^\dag & = & b^\dag a^\dag \label{ABDAG} \\
		   \e(a^\dag) & = & \overline{\e(a)} \label{EBAR} \\
		   \D(a^\dag) & = & \D(a)^\dag,\quad\forall a,b \in H \label{DDAG}  
	\end{eqnarray}
	where $\dag$ extends to a conjugation operation on all of $H \ot H$ in a natural
	way so that
	\begin{eqnarray*}
		(a \ot b)^\dag = a^\dag \ot b^\dag, \quad \forall a,b \in H.
	\end{eqnarray*}
\end{Definition}

Equations~(\ref{ADAG}) and (\ref{ABDAG}) are equivalent to the usual
definition of a conjugation operation (also referred to as a $*$-operation) on
the algebra $H$, whilst equations~(\ref{EBAR}) and (\ref{DDAG}) are the compatibility conditions
with the coalgebra structure: i.e. they determine $*$-algebra homomorphisms. 
In (\ref{EBAR}) the overbar denotes complex conjugation over $\mathbb{C}$: we adopt this
convention throughout.

A $*$-bi-algebra $H$ which admits an antipode $S$ is called a Hopf $*$-algebra. For a Hopf $*$-algebra we necessarily have for the antipode $S$~\cite{Sch93, Wor87}
\begin{Lemma}\label{SDAGGER}
        \begin{eqnarray*}
			S(a)^\dag = S^{-1}(a^\dag).
		\end{eqnarray*}
\end{Lemma}
\begin{proof}
        This follows from the uniqueness of the antipode $S$ (as the inverse of
        the identity map on $H$ under the convolution product) i.e. $S: H \rightarrow H$
        is uniquely defined by 
        \begin{eqnarray*}
			S(a_{(1)})a_{(2)}= a_{(1)}S(a_{(2)})=\e(a), \quad \forall a \in H.
		\end{eqnarray*}
        Now define $\tilde{S}:H \rightarrow H$ by
        \begin{eqnarray*}
			\tilde{S}(a)=[S^{-1}(a^\dag)]^\dag \equiv S^{-1}(a^\dag)^\dag .
		\end{eqnarray*}
        Then, since $\dag$ is compatible with $\D$,
		\begin{eqnarray*}
        	\tilde{S}(a_{(1)})a_{(2)} 
        	                          & = & \{S^{-1}[a_{(1)}^{~~\dag} S(a_{(2)}^{~~\dag})]\}^\dag \\
        	                          & = & [\overline{\e(a)}]^\dag=\e(a)
		\end{eqnarray*}		
        and similarly
        \begin{eqnarray*}
			a_{(1)} \tilde{S} (a_{(2)})= \e(a),\quad \forall a \in H.
		\end{eqnarray*}
        Thus by the uniqueness of $S$,~$\tilde{S}=S$ which is sufficient to prove the result.
\end{proof}

\begin{remark*}
	The square of the antipode, $S^2$ determines an algebra homomorphism, 
	in fact an algebra automorphism, which from lemma~\ref{SDAGGER} satisfies
	$S^2(a^\dag)=S^{-2}(a)^\dag$. Thus $S^2$ does not determine a $*$-algebra homomorphism.
\end{remark*}

For the QHA case the situation with the antipode $S$ is more complicated in view of the fact, that the antipode $S$
is no longer unique~\cite{Dri90}. Thus Lemma~\ref{SDAGGER} does not hold for QHA.

In order to formulate a suitable definition for $*$-QHA we investigate how twisting alters the $*$-structure of a Hopf $*$-algebra.
Let $H$ be a Hopf algebra and $F \in H \ot H$ an arbitrary twist. The twisted structure induced by $F$ on $H$ is no longer
a Hopf algebra but is instead a QHA. The twisted structure is obtained by setting $\Phi = 1 \otimes 1 \otimes 1$,
$\a=\b=1$ into equations (\ref{TWDEL} - \ref{TWCAN}), giving 
\begin{eqnarray}\label{HF}
	\D_F(a) & = & F \D(a)F^{-1} \no \\
     \Phi_F & = & (F \ot 1) \cdot (\D \ot 1)F \cdot (1 \ot \D) F^{-1} \cdot (1 \ot F^{-1}) \no \\
       \a_F & = & m \cdot (S \ot 1)F^{-1}, \quad \b_F=m \cdot ( 1 \ot S)F. 
\end{eqnarray}
The counit $\e$ and antipode $S$ are unchanged.

Note that $\tilde{\D}_F $ defined by
\begin{eqnarray*}
	\tilde{\D}_F(a) & = & \D_F(a^\dag)^\dag
\end{eqnarray*}
determines another co-product on $H$. It shall be shown below for the general case, that $H$
is in fact a QHA with the above co-product, with co-associator 
$\tilde{\Phi}_F=(\Phi_F^\dag)^{-1}$ and canonical elements~$\tilde{\a} =S^{-1}(\b_F)^\dag,~\tilde{\b} =S^{-1}(\a_F)^\dag $.

Since the co-product $\D$ is compatible with $\dag$ we have
\begin{eqnarray*}
	\D_F(a)^\dag & = & [F \D (a) F^{-1}]^\dag \\
				 & = & (F^\dag)^{-1} \D(a)^\dag F^\dag \\
	             & = & (F^\dag)^{-1} \D(a^\dag) F^\dag \\
	             & = & \O \D_F(a^\dag) \O^{-1} 
\end{eqnarray*}
where $\O= (FF^\dag)^{-1}$ is a self-adjoint twist. Thus
\begin{eqnarray}
	\tilde{\D}_F(a) & = & \D_F(a^\dag)^\dag = \O \D_F(a) \O^{-1},\quad\forall a \in H \no
\end{eqnarray}
so $\tilde{\D}_F$ is obtained from $\D_F$ by twisting with a (self-adjoint) twist
$\O$, or equivalently $\D_F(a)^\dag =\O \D_F(a^\dag) \O^{-1}$ as above.

Similarly for the co-associator 
\begin{eqnarray*}
	\Phi_F & = & (F \ot 1) \ (\D \ot 1)F \ (1 \ot \D) F^{-1} \ (1 \ot F^{-1}),
\end{eqnarray*}
since $\D$ is compatible with $\dag$
\begin{eqnarray*}
	(\Phi_F^\dag)^{-1} & = & ({F^\dag}^{-1} \ot 1) \ (\D \ot 1) {F^\dag}^{-1}
	\ (1 \ot \D) F^\dag \ (1 \ot F^\dag) \\
	& = &({F^\dag}^{-1} \ot 1) \ (\D \ot 1){F^\dag}^{-1}\ [(\D \ot 1) F^{-1}\ (F^{-1} \ot 1) \\
	&&  \ \Phi_F \
	(1 \ot F) \ (1 \ot \D)F] \ (1 \ot \D){F^\dag} \ ( 1 \ot {F^\dag}) \\
	& = & (\O \ot 1) \ (\D_F \ot 1) \O \ \Phi_F \ (1 \ot \D_F) \O^{-1} \
	(1 \ot \O^{-1}) \\
	& = & (\Phi_F)_\O
\end{eqnarray*}
so that $(\Phi_F^\dag)^{-1}$ is also obtained from $\Phi_F$ by twisting with the
(self-adjoint) twist $\O = (FF^\dag)^{-1}$ as above.

We are now in a position to introduce the primary object of our investigation.


\section{$*$-Quasi-Hopf algebras}
\begin{Definition}
        A QHA $H$ is called a $*$-QHA if it admits a conjugation operation $\dag$ and a twist
        $\O \in H \ot H$, called the $*$-canonical element, satisfying
        \begin{align} 
                & \e(a^\dag) =  \overline{\e(a)},\quad\forall a \in H \label{EADAG}\\
                & \D(a)^\dag =  \O \D(a^\dag) \O^{-1},\quad\forall a \in H \label{DADAG}\\
                & (\Phi^\dag)^{-1} = \Phi_\O = 
                (\O \ot 1) \ (\D \ot 1)\O \ \Phi \ (1 \ot \D) \O^{-1}
                \ (1 \ot \O^{-1}).  \label {PHIDAG}
        \end{align}
\end{Definition}
Our definition is motivated by the observation that any QHA obtained by twisting from a
Hopf $*$-algebra is a quasi-Hopf $*$-algebra ($*$-QHA). 

Following the previous section we define a new co-product $\tilde{\D}$ on $H$
by
\begin{equation}
        \tilde{\D}(a)=\D(a^\dag)^\dag,\quad \forall a \in H. \label{E28}
\end{equation}
With this co-product $H$ also determines a $*$-QHA, as will be seen below.
In view of the previous section we may have imposed the extra conditions
$\O^\dag=\O$ (i.e. $\O$ is self-adjoint) and $S(a)^\dag=S^{-1}(a^\dag),\forall a \in H $
but we will not do this below. However, we define 
\begin{Definition}
If $\O=\O^\dag$ we call a $*$-QHA $H$ self-conjugate. If the antipode $S$
satisfies 
	\begin{eqnarray*}
		S(a)^\dag=S^{-1}(a^\dag), \quad \forall a \in H
	\end{eqnarray*}
we say that $S$ is $*$-compatible.
\end{Definition}
\noindent
In general for a $*$-QHA, the antipode $S$ is not $*$-compatible, however we shall see that
$S$ is almost $*$-compatible.

Equations~(\ref{DADAG},\ref{PHIDAG}) impose strong conditions on the $*$-canonical element $\O$. Indeed,
\begin{eqnarray}
         \D(a)&=&[\D(a)^\dag]^\dag   \no \\
        & \stackrel{(\ref{DADAG})}{=}&[\O \D(a^\dag) \O^{-1}]^ \dag  \no \\
        & =&(\O^{-1})^\dag \D(a^\dag)^\dag \O^\dag \no \\
        & \stackrel{(\ref{DADAG})}{=}&(\O^{-1})^\dag 
        \O \D(a) \O^{-1} \O^\dag, \quad\forall a \in H \label{E29}
\end{eqnarray}
so that $\O^{-1}\O^\dag$, and its inverse, must commute with the co-product $\D$.
We say that $\O$ is {\it quasi-self adjoint}.

We thus have, from equation~(\ref{DADAG}),
\begin{equation}
\tilde{\D}(a)=\D(a^\dag)^\dag= \O \D(a) \O^{-1}\stackrel{(\ref{E29})}{=}
\O^\dag \D(a)(\O^\dag)^{-1} \tag{\ref{E29}$'$}
\end{equation}
or equivalently
\begin{eqnarray*}
\D(a)^\dag=\D_\O(a^\dag)=\D_{\O^\dag}(a^\dag), \quad \forall a \in H. 
\end{eqnarray*}
Thus we might expect that $\O^\dag$ is also a $*$-canonical element for $H$. This is indeed the case.
\begin{Proposition} \label{P10}
 $\O^\dag$ is also a $*$-canonical element for $H$ called the conjugate $*$-canonical element.
\end{Proposition}
\begin{proof}
        It remains to check (\ref{PHIDAG}). To this end we have
        $$(\Phi^\dag)^{-1} = \Phi_\O = (\O \ot 1) 
        \ (\D \ot 1)\O \ \Phi \ (1 \ot \D) \O^{-1} \ (1 \ot \O^{-1})
        $$
        Taking the conjugate inverse of this equation (i.e. apply $\dag$ followed
        by the inverse) and noting from (\ref{E29}$'$) 
        that $\dag \cdot \D= \tilde{\D} \cdot \dag$ ~gives
	\begin{eqnarray*}
        \Phi & = & ({\O^\dag}^{-1} \ot 1)(\tilde{\D} \ot 1){\O^\dag}^{-1} \ {\Phi^\dag}^{-1} \ (1 \ot \tilde{\D}) \O^\dag (1 \ot \O^\dag \\ 
	     & \stackrel{(\ref{E29}')} {=} & (\D \ot 1) {\O^\dag}^{-1} ({\O^\dag}^{-1} \ot 1) \ {\Phi^\dag}^{-1} (1 \ot \O^\dag) (1 \ot \D)
	 \O^\dag   
	\end{eqnarray*}
        and hence
        $$(\Phi^\dag)^{-1}=({\O^\dag} \ot 1) 
        (\D \ot 1){\O^\dag} \ \Phi
        \ (1 \ot \D) {\O^\dag}^{-1} (1 \ot {\O^\dag}^{-1})=\Phi_{\O^\dag}
        $$
        which proves the result.
\end{proof}     
It follows that for a $*$-QHA $H$ the $*$-canonical element $\O$ is not in general unique.

We now demonstrate that with the co-product $\tilde{\D}$ of equation~(\ref{E28}),
$H$ is also a $*$-QHA. In fact we have
\begin{Proposition} \label{P11}
Suppose $H$ is any QHA admitting a conjugation operation $\dag : H \rightarrow H$ satisfying only 
eq. (\ref{EADAG}). Then $H$ is a QHA with the same co-unit $\e$ but with co-product
$\tilde{\D}$, co-associator $\tilde{\Phi}=(\Phi^\dag)^{-1}$, canonical elements
$\tilde{\a}=S^{-1}(\b)^\dag, \tilde{\b}=S^{-1}(\a)^\dag$ and antipode $\tilde{S}$
defined by
\begin{equation}
\tilde{S}(a)=(S^{-1}(a^\dag))^\dag,\quad\forall a \in H. \label{E30}
\end{equation}
Moreover, if $H$ is a $*$-QHA then $H$ is also a $*$-QHA with this structure but with
canonical element $\tilde{\O}=\O^{-1}.$
\end{Proposition}
\begin{proof}
        First it is obvious that $\tilde{\D},~\e$ determine a coalgebra structure on $H$
        and are algebra homomorphisms. As to the co-associator $\tilde\Phi$ we have by applying $\dag$
        to \reff{q-co},
        $$(1 \ot \tilde{\D})\tilde{\D}(a^\dag)= (\Phi^\dag)^{-1}(\tilde{\D} \ot 1)
        \tilde{\D}(a^\dag)\Phi^\dag,\quad\forall a \in H
        $$
        which proves \reff{q-co}. As to property \reff{pentagon}, taking the conjugate
        inverse of \reff{pentagon} gives immediately
        $$(\tilde{\D} \ot 1 \ot 1)(\Phi^\dag)^{-1} \ (1 \ot 1 \ot \tilde{\D})
        (\Phi^\dag)^{-1}=({\Phi^\dag}^{-1} \ot 1) \ (1 \ot \tilde{\D} \ot 1)(\Phi^\dag)^{-1}
        \ (1 \ot {\Phi^\dag}^{-1})
        $$
        as required. Property~\reff{epsphi} is obvious, so it remains to consider
        \reff{Sphi} and \reff{Sab}. As to the former, we set
        $$ \tilde{\Phi}= \tilde{X}_\nu \ot \tilde{Y}_\nu \ot \tilde{Z}_\nu =
        (\Phi^\dag)^{-1}=  \bar{X}_\nu^\dag \ot \bar{Y}_\nu^\dag
        \ot \bar{Z}_\nu^\dag
        $$
        which implies
        \begin{eqnarray}
         \S(\tilde{X}_\nu) \tilde{\a} \tilde{Y}_\nu \tilde{\b} \S(\tilde{Z}_\nu)
        &=&\{ S^{-1}[ \bar{X}_\nu \b S(\bar{Y}_\nu) \a \bar{Z}_\nu]\}^\dag 
        \stackrel{\reff{Sphi}}{=}1 \no
        \end{eqnarray}
        and similarly setting
        \begin{eqnarray}
        \tilde{\Phi}^{-1}&=& \bar{\tilde{X}}_\nu \ot \bar{\tilde{Y}}_\nu
        \ot \bar{\tilde{Z}}_\nu =
        \Phi^\dag=  X_\nu^\dag \ot Y_\nu^\dag \ot Z_\nu^\dag \no
        \end{eqnarray}
        we have
        \begin{eqnarray}
        \bar{\tilde{X}}_\nu \tilde{\b} \S(\bar{\tilde{Y}}_\nu) \tilde{\a}
        \bar{\tilde{Z}}_\nu 
        &=&\{ S^{-1}[ S(X_\nu) \a Y_\nu \b S(Z_\nu)]\}^\dag 
        \stackrel{\reff{Sphi}}{=}1. \no
        \end{eqnarray}
        As to property \reff{Sab} we have 
        \begin{eqnarray*}
		\tilde{\D}(a) = \D(a^\dag)^\dag = (a^\dag_{~(1)})^\dag \ot (a^\dag_{~(2)})^\dag
		\end{eqnarray*}
        so that
        \begin{eqnarray}
         \S[(a^\dag_{~(1)})^\dag] \tilde{\a}(a^\dag_{~(2)})^\dag 
        &=& S^{-1}[(a^\dag_{~(1)}) \b S(a^\dag_{~(2)})]^\dag \no \\
        &=& \e(a)S^{-1}(\b)^\dag=\e(a) \tilde{\a} \no
        \end{eqnarray}
        and similarly for $\tilde{\b}$ as required. This proves that $H$ gives rise
        to a QHA under the given structure.

        Finally if $H$ is a $*$-QHA with $*$-canonical element $\O$ then $H$ is also a 
        $*$-QHA under the above structure but with $*$-canonical element $\O^{-1}$. To see this
        we have from equation~(\ref{E29}$'$),
        $$\tilde{\D}(a)^\dag=\D(a^\dag)= \O^{-1}\D(a)^\dag \O = \O^{-1} \tilde{\D}(a^\dag)
        \O,\quad\forall a \in H$$
        which proves (\ref{DADAG}), while for (\ref{PHIDAG}) we have
        $$\tilde{\Phi}_{\O^{-1}}=(\Phi^\dag)^{-1}_{\O^{-1}}=(\Phi_{\O})_{\O^{-1}}=\Phi_{\O^{-1} \O}=\Phi 
        = {(\tilde{\Phi})^\dag}^{-1}$$
        so that $\O^{-1}$ is a $*$-canonical element for this structure thus making it a 
        $*$-QHA.
\end{proof}

When $H$ admits a conjugation operation $\dag$ satisfying (\ref{EADAG})~
it ensures that with the structure of proposition \ref{P11}, $H$ is also a QHA.
Conditions (\ref{DADAG}, \ref{PHIDAG}) are equivalent to this
QHA structure being obtainable, up to equivalence modulo $(S, \a,\b)$, by twisting with
$\O$.

We now demonstrate that the category of $*$-QHAs is invariant under
twisting, as is the sub-category of self-conjugate $*$-QHAs.
This latter observation is important as it demonstrates that we cannot obtain
a self-conjugate $*$-QHA from a non-self-conjugate one by twisting.

\begin{Theorem} Let $H$ be a (self-conjugate) $*$-QHA with $*$-canonical element
$\O$ and $F \in H \ot H$ an arbitrary twist. Then $H$ is also a (self-conjugate)
$*$-QHA with the twisted structure of equations (\ref{TWDEL}) with $*$-canonical element 
$\O_F=(F^\dag)^{-1} \O F^{-1}$. Moreover if the antipode $S$ is $*$-compatible
then it is $*$-compatible under this twisted structure. \label{T5}
\end{Theorem}
\begin{proof}
        It suffices to prove (\ref{DADAG}, \ref{PHIDAG}). For the twisted
        co-product we have
        \begin{eqnarray}
        \D_F(a)^\dag&=&[F \D(a)F^{-1}]^\dag= (F^{-1})^\dag \D(a)^\dag F^\dag \no \\
        &\stackrel{(\ref{DADAG})}{=}& 
        (F^{-1})^\dag \O \D(a^\dag) \O^{-1} F^\dag \no \\
        &=& \O_F \D_F(a^\dag) \O_F^{-1} \no
        \end{eqnarray}
        with $\O_F= (F^{-1})^\dag \O F^{-1}$ as stated. As to the co-associator we have
        $$\Phi_F=(F \ot 1) \ (\D \ot 1)F \ \Phi
        \ (1 \ot \D) F^{-1} \ (1 \ot F^{-1})$$
        so that taking the conjugate inverse gives
        $$(\Phi_F^\dag)^{-1} =({F^\dag}^{-1} \ot 1) \ (\tilde{\D} \ot 1){F^\dag}^{-1}
        \ {\Phi^\dag}^{-1}
        \ (1 \ot \tilde{\D}) F^\dag \ (1 \ot F^\dag)$$
        with $\tilde{\D}$ as in equation (\ref{E28})
        [also cf equation (\ref{E29}$~'$)] and where
        $$(\Phi^\dag)^{-1}\stackrel{(\ref{DADAG})}{=} 
        \Phi_\O=(\O \ot 1) \ (\D \ot 1)\O \ \Phi
        \ (1 \ot \D) \O^{-1} \ (1 \ot \O^{-1}).$$
        Thus by equation (\ref{E29}$'$)
        \begin{eqnarray}
        (\Phi_F^\dag)^{-1}&=&({F^\dag}^{-1}\ot 1) \ (\D_\O \ot 1){F^\dag}^{-1} \ (\O \ot 1)
        \ (\D \ot 1) \O  \no \\
        &&\ \Phi \ (1 \ot \D) \O^{-1} \ (1 \ot \O^{-1}) \ (1 \ot \D_\O)F^\dag
        \ (1 \ot F^\dag)  \no \\
        &=&({F^\dag}^{-1} \O \ot 1) \ (\D \ot 1)({F^\dag}^{-1} \O) \no \\
        &&\ \Phi
        \ (1 \ot \D)(\O^{-1}F^\dag) \ (1 \ot \O^{-1}F^\dag) \label{E31} \\
        &=&({F^\dag}^{-1} \O \ot 1) \ (\D \ot 1)({F^\dag}^{-1} \O) \
        (\D \ot 1) F^{-1} \ (F^{-1} \ot 1) \no \\
        &&\ \Phi_F \ (1 \ot F) \ (1 \ot \D) F
        \ (1 \ot \D)(\O^{-1}F^\dag) \ (1 \ot \O^{-1}F^\dag) \no \\
        &=&(\O_F \ot 1) \ (\D_F \ot 1)\O_F \ \Phi_F \ (1 \ot \D_F){\O_F}^{-1}
        \ (1 \ot {\O_F}^{-1}) \no \\
        &=&(\Phi_F)_{\O_F} \no
        \end{eqnarray}
        with $\O_F=(F^\dag)^{-1} \O F^{-1}$ as required. Thus under the twisted structure
        induced by $F$, $H$ is a $*$-QHA with $*$-canonical element $\O_F$ as stated.
        If moreover $H$ is  self conjugate, so that $\O$ is self adjoint, so too
        is $\O_F$ which implies $H$ is also a self-conjugate $*$-QHA under the twisted
        structure. Finally $*$-compatibility of the antipode $S$ is obviously twist
        invariant since $S$ remains unchanged under twisting.
\end{proof}
We refer to the twisted structure above as the twisted $*$-QHA induced by $F$. The above result has a number of
interesting consequences to which we now turn.

\begin{Proposition} \label{P12}
Let $H$ be a $*$-QHA with $*$-canonical element $\O$. Then $H$ is also a $*$-QHA under the 
opposite structure of proposition \ref{P1} with $*$-canonical element $\O^T=T \cdot \O$.
\end{Proposition}
\begin{proof}
        Recall that $H$ is a QHA under the opposite structure with co-product $\D^T=T \cdot \D$,
        co-associator $\Phi^T=\Phi^{-1}_{321}$ and antipode $S^{-1}$, with the same co-unit.
        To prove this gives rise to a $*$-QHA it suffices to prove (\ref{DADAG}, \ref{PHIDAG}).
        For the co-product we have,
        $$\D^T(a)^\dag=T \cdot [\D(a)^\dag]=T \cdot[\O \D(a^\dag)\O^{-1}]
        =\O^T \D^T(a^\dag)(\O^T)^{-1},\quad\forall a \in H$$
        as required. For the co-associator we have
        \begin{align}
        {(\Phi^T)^\dag}^{-1}=\Phi^\dag_{321}. \label{S8}
        \end{align}
        Now since
        $$(\Phi^\dag)^{-1}=\Phi_\O=
        (\O \ot 1) \ (\D \ot 1)\O \ \Phi \ (1 \ot \D) \O^{-1} \ (1 \ot \O^{-1})$$
        we have
        $$\Phi^\dag=
        (1 \ot \O) \ (1 \ot \D)\O \ \Phi^{-1} \ (\D \ot 1) \O^{-1} \ (\O^{-1} \ot 1)$$
        and hence
        $${(\Phi^T)^\dag}^{-1}\stackrel{(\ref{S8})}{=}
        [(1 \ot \O) \ (1 \ot \D)\O \ \Phi^{-1} \ (\D \ot 1) \O^{-1}
        \ (\O^{-1} \ot 1)]_{321}$$
        $$=(\O^T \ot 1) \ (\D^T \ot 1 )\O^T \ \Phi^T \ (1 \ot \D^T) {\O^T}^{-1}
        \ (1 \ot {\O^T}^{-1} )= (\Phi^T)_{\O^T}$$
        so ${(\Phi^T)^\dag}^{-1}$ is obtained from $\Phi^T$ by twisting with $\O^T$
        under the opposite structure. Thus $H$ is also a $*$-QHA with $*$-canonical
        element $\O^T$ under the opposite structure as required.
        If moreover $H$ is self-conjugate, so that $\O$ is self-adjoint, so too
        is $\O^T$. Thus under the opposite structure, a self-conjugate $*$-QHA, is
        also self-conjugate. Obviously if the antipode $S$ of $H$ is $*$-compatible
        so too is the antipode $S^{-1}$ for the opposite structure.
\end{proof}

We have already seen that the $*$-canonical element $\O$ for a $*$-QHA is not unique,
since $\O^\dag$ also gives rise to a $*$-canonical element . We thus conclude this
section with the following observation on the uniqueness, and existence of
$*$-canonical elements.

Let $F,G \in H \otimes H$ be twists on $H$. The composite twist $FG$ is given by first twisting with $G$ and then
twisting $H_G$ by $F$ so that
\begin{eqnarray*}
X_{FG} = (X_G)_F
\end{eqnarray*}
where $X$ is one of $\D,\Phi,\a,\b,\R$. A twist $C \in H \ot H$ which preserve the QBA structure on $H$, so that
\begin{eqnarray*}
\quad \D_C(a) &=& \D(a),\quad \forall a \in H \\
\quad  \Phi_C &=& \Phi \label{QCC1}
\end{eqnarray*}
is called a {\em compatible} twist~\cite{mdgtl}. The set of compatible twists is a subgroup of the group of all twists
on $H$.

\begin{Theorem} \label{T6}
Let $H$ be a $*$-QHA with $*$-canonical element $\O$. Then $\Gamma \in H \ot H$ is also a $*$-canonical 
element for $H$ if and only if there exists a (unique) compatible twist $F \in H \ot H$
such that $\Gamma =\O F$.
\end{Theorem}
\begin{proof}
Follows from a direct computation using the composition laws for twists.
\end{proof}
\begin{Corollary*}
For a $*$-QHA $H$, there is a one to one correspondence between $*$-canonical elements and
compatible twists on $H$.
\end{Corollary*}

In particular there must exist a compatible twist
$C \in H \ot H$ such that $\O^\dag=\O C$. Thus we see that $\O$ is almost self-adjoint,
hence the term quasi-self adjoint.
As will be seen below the explicit choice of $*$-canonical element $\O$ has no effect on
the algebraic properties of $*$-QHAs, due to the special nature of compatible twists.

The existence of a conjugation operation on a $*$-QHA and the properties (\ref{EADAG}-
\ref{PHIDAG})
imply some interrelationships between $\dag$ and the algebraic structure of $H$ to 
which we now turn.


\section{Compatibility of $*$ and algebra properties}

A QHA differs from a Hopf algebra in that the antipode $S$ and its corresponding canonical elements $\a,\b$
are not unique. Nevertheless, the antipode and its corresponding canonical elements are almost unique
as the following result due to Drinfeld~\cite{Dri90} shows.

\begin{Theorem}\label{T1} Suppose $H$ is also a QHA with antipode $\S$ and canonical elements $\tilde \a ,\tilde \b$.
Then there exists a unique invertible $v \in H$ such that
\begin{equation*}
v\a=\tilde \a~~,\tilde \b v=\b~~,\S (a)=v S(a)v^{-1}~~,\forall a\in H.
\end{equation*}
Explicitly
\begin{eqnarray*}
(i)& v= \S(X_\nu)\tilde\a Y_\nu\b S{(Z_\nu)}=\S
{(S^{-1}(\X_\nu))}\S {(S^{-1}(\b))}\S{(\Y_\nu)}\tilde\a\Z_\nu\\
(ii)& v^{-1}= S{(X_\nu)}\a Y_\nu\tilde\b\S {(Z_\nu)}=
\X_\nu\tilde\b\S{(\Y_\nu)}\S{(S^{-1}(\a))} \S
{(S^{-1}{(\Z_\nu)})} \\
\end{eqnarray*}
\end{Theorem}
For arbitrary invertible $v \in H$, the triple $(\S, \tilde\a, \tilde\b)$ defined by
\begin{eqnarray*}
\S(a) = vS(a)v^{-1}, \quad \tilde\a = v\a, \quad \tilde\b = \beta v^{-1}
\end{eqnarray*}
satisfies equations ~\reff{Sphi},~\reff{Sab} and hence gives rise to an antipode $\S$ with corresponding
canonical elements $\tilde\a,\tilde\b$. There is thus a $1-1$ correspondence between
triples $(\S,\tilde \a,\tilde \b)$ and invertible $v \in H$. We say that these structures are equivalent (modulo $(S,\a,\b)$)
as they give rise to equivalent QHA structures.

Proposition \ref{P11} shows that $H$ is a $*$-QHA with co-unit $\e$, co-product $\tilde{\D}=
\D_\O=\D_{\O^\dag}$, co-associator $(\Phi^\dag)^{-1}=\Phi_\O=\Phi_{\O^\dag}$, canonical
elements $\tilde{\a}=S^{-1}(\b)^\dag,~\tilde{\b}=S^{-1}(\a)^\dag$ and with antipode
$\S$ given by equation (\ref{E30}). On the other hand, from equations (\ref{EADAG}-\ref{PHIDAG}),
$H$ is also a $*$-QHA under the twisted structure induced by $\O$ (or
$\O^\dag$) with the same co-unit, co-product and co-associator but with antipode $S$ and
twisted canonical elements given by equation (\ref{TWCAN}) 
$$\a_\O=m \cdot ( 1 \ot \a)(S \ot 1) \O^{-1},~ \b_\O= m \cdot (1 \ot \b)(1 \ot S)\O$$
or
$$\a_{\O^\dag}=m \cdot ( 1 \ot \a)(S \ot 1) (\O^\dag)^{-1},~ \b_{\O^\dag}= 
m \cdot (1 \ot \b)(1 \ot S)\O^\dag.$$
Hence these structures must be equivalent. We have immediately from Theorem \ref{T1}
\begin{Proposition} \label{P14}
There exists a unique invertible $w \in H$ such that
\begin{eqnarray} \label{E32}
&(i)&wS^{-1}(\b)^\dag=\a_\O,~\b_\O w=S^{-1}(\a)^\dag \no \\
&(ii)&S(a)=w \S(a) w^{-1},\quad\forall a \in H. 
\end{eqnarray}
Explicitly
\begin{eqnarray*}
w&=& S(\bar{X}_\nu^\dag) \a_\O \bar{Y}_\nu^\dag S^{-1}(\a)^\dag \S(\bar{Z}_\nu^\dag) \\
&=& S(\S^{-1}({X}_\nu^\dag)) S(\S^{-1}[S^{-1}(\a)^\dag]) S({Y}_\nu^\dag) 
\a_\O {Z}_\nu^\dag \\
&& \\
w^{-1}&=& \S(\bar{X}_\nu^\dag) S^{-1}(\b)^\dag \bar{Y}_\nu^\dag \b_\O S(\bar{Z}_\nu^\dag) \\
&=& {X}_\nu^\dag \b_\O  S({Y}_\nu^\dag) S(\S^{-1}[S^{-1}(\b)^\dag])
S(\S^{-1} ({Z}_\nu^\dag)).
\end{eqnarray*}
\end{Proposition}

Above we used the fact that the co-associator for the QHA we are considering is 
$(\Phi^\dag)^{-1}$ together with the antipode $\S$ and canonical elements $S^{-1}(\b)^\dag,
S^{-1}(\a)^\dag$ respectively. We then applied Theorem \ref{T1} to this structure with
$(\S, \tilde\a,\tilde\b) \equiv (S,\a_\O,\b_\O)$.

\begin{Corollary}\label{C14_1}
$$w^\dag =  S^{-1}(\bar{Z}_\nu) S^{-1}(\a) \bar{Y}_\nu \a_{\O}^\dag \S^{-1}(\bar{X}_\nu)$$ 
$$(w^{-1})^\dag =  \S^{-1}(\bar{Z}_\nu) \b_\O^\dag \bar{Y}_\nu S^{-1}(\b) S^{-1}(\bar{X}_\nu).$$ 
\end{Corollary}

\begin{Corollary} \label{C14_2}
$S$ is $*$-compatible i.e. $S(a)^\dag=S^{-1}(a^\dag)$ or equivalently  $\S(a) \equiv
[S^{-1}(a^\dag)]^\dag$ $=S(a),$ $~\forall a \in H$ if and only if $w$ as above is a central element.
\end{Corollary}
\setcounter{Corollary}{0}
The results above, particularly equation~(\ref{E32})(ii) and Corollary~\ref{C14_2}
might be thought to depend on the $*$-canonical element $\O$. To see this is not the case,
let $\Gamma$ be another $*$-canonical element so $\Gamma= \O C$, for some compatible twist $C \in H \ot H$. 
The corresponding twisted canonical elements are
$$\a_\Gamma=\a_{\O C} = (\a_C)_\O,~~\b_\Gamma=\b_{\O C}=(\b_C)_\O~.$$
From Theorem~\ref{T1} there exists a unique invertible element $z \in H$ such that
$$\a_C=z \a, \quad \b_C=z^{-1} \b$$
with
$$ S(a) = zS(a)z^{-1},\quad \forall a \in H.$$
The element $z$ is thus central.
Now,
$$\a_\Gamma = z \a_\O,~~\b_\Gamma=z^{-1} \b_\O.$$
The corresponding $w$-operator, given by replacing $\a_\O,\b_\O$ with $\a_\Gamma,\b_\Gamma$
respectively, is thus given by
$$w_\Gamma=zw,~~w_\Gamma^{-1}=z^{-1}w^{-1}$$
so that, in particular
$$S(a)=w\S(a)w^{-1}=w_\Gamma \S(a) w_\Gamma^{-1} $$
the latter equality holding identically. Thus the results of Proposition~\ref{P14} and
its corollaries are independent (modulo an invertible central element) of the canonical element chosen.
\begin{flushright}$\Box$\end{flushright}

We have shown previously~\cite{mdgtl} that the $v$ operator of Theorem~\ref{T1} is
universal i.e. unchanged under twisting by an arbitrary twist $F$, so that for
any operator $v$ arising from the application of Theorem~\ref{T1} we have $v_F=v$. Since
the $w$ operator arises precisely in this way, it follows that
\begin{Theorem}
The operator $w$ is universal, i.e. twist invariant.
\end{Theorem}

\begin{remark*}
The universality of $w$ can also be shown by direct calculation. The operator $w$ can be expressed in the following simpler form,
\begin{equation}\label{altw}
w= S(X_\nu)\a Y_\nu S^{-1}(\a_{\O^\dag})^\dag \S(Z_\nu).
\end{equation}
\end{remark*}

The results of Proposition~\ref{P14} have a number of interesting consequences which
we summarise below:
\begin{Lemma} \label{L7}(Notation as above)
\begin{eqnarray}
&(i)&\S(S^{-1}(a))=w^{-1}aw,\quad\forall a \in H \no \\
&(ii)&S(\S^{-1}(a))=waw^{-1},\quad\forall a \in H \no   \\
&(iii)&S(\S^{-1}(w))=\S(S^{-1}(w))=w \no
\end{eqnarray}
so that
\begin{eqnarray}
&&S^{-1}(w)^\dag=S(w^\dag)~~ {\rm or} ~~S^{-1}(w)=\S^{-1}(w) \label{E33} 
\end{eqnarray}
and similarly for $w^{-1}$.
\begin{eqnarray}
&(iv)&S^{-1}(a)=S^{-1}(w)\S^{-1}(a)S^{-1}(w^{-1}) \no \\
&&~~~~~~~~=\S^{-1}(w)\S^{-1}(a)\S^{-1}(w^{-1}),~~~\forall a \in H. \no
\end{eqnarray}
If z $\in H$ is a central element then 
\begin{eqnarray}
&(v)&S^{-1}(z)=\S^{-1}(z) \no
\end{eqnarray}
so that 
\begin{eqnarray}
&&S^{-1}(z)^\dag=S(z^\dag)~. \no
\end{eqnarray}
\end{Lemma}
\begin{proof}
        (i), (ii) and (v) follow directly from equation (\ref{E32})(ii). Part
	(iii) is a direct consequence of (i) and (ii). Applying $S^{-1}$ to
	part (ii) gives (iv).
\end{proof}

Since $\O^\dag$ is also a $*$-canonical element for $H$ we may replace $\O$ with $\O^\dag$
in proposition \ref{P14} to give 
\vskip 3mm
\noindent
{\bf Proposition \ref{P14}$'$~~} 
{\it There exists a unique invertible $\bar w \in H$ such that
\begin{eqnarray}
&(i)&\bar w S^{-1}(\b)^\dag= \a_{\O^\dag},~~\b_{\O^\dag} \bar w = S^{-1} ( \a)^\dag \no \\
&(ii)&S(a)= \bar w \S (a) \bar w^{-1},\quad\forall a \in H. \no
\end{eqnarray}
Explicitly $\bar w$ is given as in proposition \ref{P14} with $\O$ replaced by $\O^\dag$.}
\begin{Corollary*}
$c=w^{-1} \bar w = \bar w w^{-1}$ is a central element with inverse 
\begin{equation}
c^{-1}= w \bar w ^{-1} = \bar w ^{-1} w. \label{E34}
\end{equation}
\end{Corollary*}

Thus the results of lemma~\ref{L7} also hold for $\bar w$. 
In view of the definition~(\ref{E30}) of $\S$; i.e.
$$\S(a)=S^{-1}(a^\dag)^\dag,\quad\forall a \in H$$
the canonical elements of proposition~\ref{P11} may be written
$$S^{-1}(\b)^\dag=\S(\b^\dag),~~S^{-1}(\a)^\dag=\S(\a^\dag)~.$$
Also, by taking the Hermitian conjugate of (\ref{E32})(ii), equation~(\ref{E32})
may be written as
\begin{align}
(i)&\quad w\S(\b^\dag)= \a_\O,~~\b_\O w=\S(\a^\dag) \no \\
(ii)&\quad S(a)^\dag=(w^\dag)^{-1}  S^{-1}(a^\dag) w^\dag,\quad\forall a \in H \tag{\ref{E32}$'$}
\label{E32'}
\end{align}
and similarly for $\O^\dag$, with $w$ replaced by $\bar w$.

It might be thought that the operators $w,\bar w$ of propositions~(\ref{P14},\ref{P14}$'$) respectively,
are directly related. This turns out to be the case. We first need
\begin{Lemma} \label{L8}
(notation as above):
\begin{eqnarray*}
&(i)& \a_\O^{\ \dag}=S^{-1}(\b)S^{-1}(\bar w), \quad  \b_\O^{\ \dag}=S^{-1}(\bar w^{-1})S^{-1}(\a)  \\
&(ii)& \a_{\O^\dag}^{\ \dag}=S^{-1}(\b)S^{-1}(w), \quad \b_{\O^\dag}^{\ \dag}=S^{-1}(w^{-1}) S^{-1}(\a). 
\end{eqnarray*}
\end{Lemma}
\begin{proof}
        By symmetry, it suffices to prove (i). Below we write (summation
        over repeated indices assumed)
        $$\O= \O_i \ot \O^i,\quad\O^{-1}= \bar{\O}_i \ot \bar{\O}^i.$$
        We have
        $$\a_\O^{\ \dag}=[S(\bar{\O}_i ) \a \bar{\O}^i]^\dag=(\bar{\O}^i)^\dag \a^\dag 
        \S^{-1}(\bar{\O}_i^\dag)$$
        where
        $$\a^\dag \stackrel{(\ref{E32}')(i)}{=}\S^{-1}(\b_{\O^\dag} \bar w)
        =\S^{-1}(\bar w)\S^{-1}(\b_{\O^\dag})
        \stackrel{(\ref{E33})(iii)}{=}S^{-1}(\bar w)\S^{-1}(\b_{\O^\dag})$$
        and we have used the fact that eq.~(\ref{E32}$'$) also holds for $\O^\dag$ with $w$
        replaced by $\bar w$. Thus
        \begin{eqnarray}
        \a_\O^{\ \dag} &=& (\bar{\O}^i)^\dag S^{-1}(\bar w)
        \S^{-1}(\b_{\O^\dag})\S^{-1}(\bar{\O}_i^\dag) \no \\
        &\stackrel{(\ref{E33})(iv)}{=}& (\bar{\O}^i)^\dag S^{-1}(\b_{\O^\dag})
        S^{-1}(\bar{\O}_i^\dag) S^{-1}(\bar w) \no \\
	&=& S^{-1}(\b)S^{-1}(\bar w).\no
        \end{eqnarray}
        The result for $\b_\O^{\ \dag}$ is proved in a similar way.
\end{proof}

We are now in a position to compute $w ^\dag$. We have
\begin{Proposition}\label{P15}
(notation as above)
$$w^\dag=S^{-1}(\bar w),\quad\bar w^\dag=S^{-1}(w).$$
In particular, for the central element $c$ of equation (\ref{E34}), we have 
$$c^\dag=S^{-1}(c^{-1}).$$
\end{Proposition}
\begin{proof}
        Using corollary \ref{C14_1} to proposition \ref{P14} we have
        \begin{eqnarray}
        w^\dag&=& S^{-1}(\bar{Z}_\nu) S^{-1}(\a) 
        \bar{Y}_\nu \a_\O^{\ \dag} \S^{-1}(\bar{X}_\nu) \no \\
        &\stackrel{lemma (\ref{L8})(i)}{=}& 
        S^{-1}(\bar{Z}_\nu)S^{-1}(\a) \bar{Y}_\nu 
        S^{-1}(\b)S^{-1}(\bar w) \S^{-1}(\bar{X}_\nu)\no \\
        &\stackrel{(\ref{E33})(iv)}{=}&
        S^{-1}(\bar{Z}_\nu)S^{-1}(\a) \bar{Y}_\nu 
        S^{-1}(\b) S^{-1}(\bar{X}_\nu)S^{-1}(\bar w)\no \\
        &=&S^{-1}[ \bar{X}_\nu \b S(\bar{Y}_\nu) \a \bar{Z}_\nu]S^{-1}(\bar w)
        \stackrel{(\ref{Sphi})}{=}S^{-1}(\bar w).\no 
        \end{eqnarray}
        Thus
        $$\bar w= S(w^\dag)$$
        which implies
        $$\bar w^\dag=S(w^\dag)^\dag=\S^{-1}(w)\stackrel{(\ref{E33})(iii)}{=}S^{-1}(w)$$
        and similarly $$w^\dag=S^{-1}(\bar{w}).$$
        Finally, as to the central element $c=w^{-1} \bar w$ of equation (\ref{E34}) we
        have, from the above
        $$c^\dag=\bar w^\dag (w^\dag)^{-1}=S^{-1}(w)S^{-1}(\bar w^{-1})
        =S^{-1}(\bar w^{-1} w)=S^{-1}(c^{-1})$$
        which proves the result.
\end{proof}
\begin{Corollary*}
$$\bar w = S(w^\dag)= S(\S^{-1}(\bar{X}_\nu)) S(\a_\O^\dag)
S(\bar{Y}_\nu) \a \bar{Z}_\nu $$   
$$\bar w^{-1} = S(w^\dag)^{-1}= \bar{X}_\nu \b S(\bar{Y}_\nu)
S(\b_\O^\dag) S(\S^{-1}(\bar{Z}_\nu)) $$  
and similarly for $w,w^{-1}$ with $\O$ replaced by $\O^\dag$.
\end{Corollary*}   
\begin{proof}
        Follows from applying the result above and $S$ to Corollary~\ref{C14_1} of
        Proposition~\ref{P14}.
\end{proof}

In the case $S$ is $*$-compatible so that $S=\S$, or equivalently $w, \bar w$
are both central, the corollary above reduces to
$$\bar w = \bar{X}_\nu S(\a_\O^\dag) S(\bar{Y}_\nu) \a \bar{Z}_\nu $$   
$$\bar w^{-1} = \bar{X}_\nu \b S(\bar{Y}_\nu) S(\b_\O^\dag) \bar{Z}_\nu $$   
and similarly for $w,w^{-1}$ with $\O$ replaced by $\O^\dag$. This gives a useful
expansion directly in terms of $\Phi^{-1}$. 

In the case the $*$-QHA $H$ is self-conjugate, so that $\O=\O^\dag$ and
$w=\bar w$, the result of proposition~\ref{P15} gives
$$w^\dag=S^{-1}(w)$$
while the central element $c$ of equation~(\ref{E34}) is obviously trivial.

We conclude this section with a simple observation concerning conjugation of the
twisted operators of equations~(\ref{TWDEL}), of use below.
\begin{Lemma}\label{L9}
Let $F \in H \ot H$ be a twist on a $*$-QHA $H$ with $*$-canonical element $\O$. Then
\begin{eqnarray}
&(i)& \D_F^\dag(a)=\D_{(F^\dag)^{-1}\O}(a^\dag),\quad\forall a \in H \no \\
&(ii)&(\Phi_F^\dag)^{-1}=\Phi_{(F^\dag)^{-1}\O},
\quad(\Phi_F^\dag)=\Phi_{(F^\dag)^{-1}\O}^{-1} \no \\
&(iii)& \a_F^{\ \dag} =S^{-1}[\b_{(F^\dag)^{-1}\O}]S^{-1}(w),
\quad \b_F^{\ \dag}=S^{-1}(w^{-1})S^{-1}[\a_{(F^\dag)^{-1}\O}]. \no
\end{eqnarray}
\end{Lemma}
\begin{proof}
        (i)
        \begin{eqnarray} \D_F(a)^\dag &=& [F\D(a)F^{-1}]^\dag=(F^{-1})^\dag \D(a)^\dag F^\dag  \no \\
        &=& (F^{-1})^\dag\O\D(a^\dag)\O^{-1}F^\dag \no \\
        &=& \D_{(F^\dag)^{-1}\O}(a^\dag),\quad\forall a \in H. \no
        \end{eqnarray}
        (ii) From equation (\ref{E31}) in the proof of Theorem \ref{T5}, we have
        \begin{align}
        (\Phi_F^\dag)^{-1}&=
        ({F^\dag}^{-1} \O \ot 1) \cdot (\D \ot 1)({F^\dag}^{-1} \O) \cdot \Phi 
        \cdot (1 \ot \D)(\O^{-1}F^\dag) \cdot (1 \ot \O^{-1}F^\dag) \no \\  
        &=\Phi_{(F^\dag)^{-1}\O}. \no
        \end{align}
        (iii) First set (summation over repeated indices)
        $$F= f_i \ot f^i,\quad F^{-1}= \bar f_i \ot \bar f^i.$$
        Then
        \begin{eqnarray}
        \a_F^{\ \dag}&=&[S(\bar f_i) \a\bar f^i]^\dag=(\bar f^i)^\dag \a^\dag \S^{-1} (\bar f_i^\dag) \no \\
        &\stackrel{(\ref{E32}')(i)}{=}&(\bar f^i)^\dag \S^{-1}(\b_\O w)\S^{-1}(\bar f_i^\dag) \no \\
        &\stackrel{(\ref{E33})(iii)}{=}&(\bar f^i)^\dag S^{-1}(w)\S^{-1}(\b_\O)\S^{-1}(\bar f_i^\dag) \no\\
        &\stackrel{(\ref{E33})(iv)}{=}&(\bar f^i)^\dag S^{-1}(\b_\O) S^{-1}(\bar f_i^\dag) S^{-1}(w) \no \\
        &=&S^{-1}[\b_{(F^\dag)^{-1} \O}] S^{-1}(w). \no
        \end{eqnarray}
        The proof for $\beta_F^{\ \dag}$ is similar.
\end{proof}

\begin{remark*}
In the case $F=\O^\dag$, the results of Lemma (\ref{L9})(iii) reduce to those of Lemma~(\ref{L8})(iii,iv).
\end{remark*}

We now turn our attention to some more advanced results on the compatibility of the
algebraic and $*$-properties of $*$-QHAs.


\section{Conjugation of the Drinfeld Twist}

Observe that $\D'$ defined by
\begin{equation}
\D'(a) = {(S\ot S)}\D^T{(S^{-1}(a))}, \quad \forall a \in H \label{E8}
\end{equation}
also determines a co-product on $H$. 
\begin{Proposition}\label{P2}
Let $H$ be a QHA, then $H$ is also a QHA with the same co-unit $\e$
and antipode $S$ but with co-product $\D'$,  co-associator $\Phi'={(S\ot S\ot S)}\Phi_{321}$
and canonical elements $\a'= S(\b),~\b'=S(\a)$. 
\end{Proposition}

Drinfeld has proved that this QHA structure is obtained by twisting with the Drinfeld twist, herein denoted $F_\d$,
given explicitly by
\begin{eqnarray*}
&(i)&F_\d={(S\ot S)}\D^T{(X_\nu)} \cdot \gamma \cdot \D {(Y_\nu\b S{(Z_\nu)})} \\
&\mbox{}& \quad =\D'{(\X_\nu \b S{(\Y_\nu)})} \cdot \gamma \cdot\D{(\Z_\nu)} 
\end{eqnarray*}
where
\begin{eqnarray}
&(ii)&\gamma= S{(B_i)}\a C_i\ot S{(A_i)}\a D_i \no
\end{eqnarray}
with
\begin{eqnarray}
&(iii)& A_i\ot B_i\ot C_i\ot D_i=\left\{
\begin{array}{l}
(\Phi^{-1}\ot 1) \cdot(\D \ot 1 \ot 1) \Phi \\
\textrm{or} \\
(1 \ot \Phi)\cdot (1 \ot 1 \ot \D)\Phi^{-1}~.
\end{array} 
\right. \no \\ \label{E9}
\end{eqnarray}
The inverse of $F_\d$ is given explicitly by
\begin{eqnarray}
&(i)& F_\d^{-1}=\D{(\X_\nu)} \cdot \bar \gamma \cdot \D'
{(S{(\Y_\nu)}\a\Z_\nu)}\no\\
&\mbox{}& \quad = \D{(S{(X_\nu)} \a Y_\nu)} \cdot \bar \gamma \cdot{(S\ot S)}\D^T
{(Z_\nu)}\no
\end{eqnarray}
where
\begin{eqnarray}
&(ii)&\bar\gamma=\bar A_i\b S{(\bar D_i)}\ot\bar B_i\b S{(\bar C_i)} \no
\end{eqnarray}
with
\begin{eqnarray}
&(iii)& \bar A_i\ot \bar B_i\ot \bar C_i\ot \bar D_i=\left\{
\begin{array}{l}
(\D \ot 1 \ot 1) \Phi^{-1} \cdot (\Phi\ot 1) \\
\textrm{or} \\
(1 \ot 1 \ot \D)\Phi\cdot (1 \ot \Phi^{-1})~. 
\end{array} 
\right. \no \\ \label{E10}
\end{eqnarray}

Replacing $S$ with $S^{-1}$ in eq.~(\ref{E8}) we obtain yet another co-product $\D_0$ on $H$:
\begin{equation}
\D_0(a)={(S^{-1}\ot S^{-1})}\D^T{(S(a))},\quad \forall a\in H. \tag{\ref{E8}$'$}
\end{equation}
We have the following analogue of proposition~\ref{P2}, the proof of which
parallels that of~\cite{oqhsa} proposition 4, but with $S$ and $S^{-1}$ interchanged:
\vskip 3mm
\noindent{\bf Proposition \ref{P2}$'$~~}
{\it $H$ is also a QHA with the same co-unit $\e$ and antipode $S$ but with co-product
$\D_0$, co-associator $\Phi_0={(S^{-1}\ot S^{-1}\ot S^{-1})}\Phi_{321}$ and 
canonical elements $\a_0=S^{-1}(\b),~\b_0=S^{-1}(\a)$ respectively.
}
\vskip 3mm
By symmetry we would expect this structure to be obtainable twisting. 
Indeed we have
\begin{Theorem}\label{T2}: The QHA structure of proposition~\ref{P2}$'$ is 
obtained by twisting with
\begin{equation}
F_0\equiv{(S^{-1}\ot S^{-1})}F_\d^T\label{E13}
\end{equation}
herein referred to as the second Drinfeld twist, where $F_\d$ is the Drinfeld twist and $F_\d^T=T\cdot F_\d.$
\end{Theorem}

Throughout we assume that $H$ is a $*$-QHA with $*$-canonical element
\begin{eqnarray*}
\O= \O_i \ot \O^i,\quad\O^{-1}= \bar \O_i \ot \bar \O^i
\end{eqnarray*}
(summation over repeated indices). In view of Theorem~\ref{T5}, $H$ is also a 
$*$-QHA under the QHA structures of propositions~\ref{P2}, \ref{P2}$'$ induced by
twisting with the Drinfeld twists $F_\delta$ and $F_0$ respectively, with $F_\delta$
as in equation (\ref{E9}) and $F_0$ as in equation (\ref{E13}). Further from 
Theorem~\ref{T5}, the $*$-canonical elements for these QHAs are given by
\begin{equation}\label{E35}
\O'=(F_\delta^\dag)^{-1} \O F_\delta^{-1},\quad\O_0=(F_0^\dag)^{-1} \O F_0^{-1}.
\end{equation}
It is one of the aims below to obtain the operators of equation~(\ref{E35}) explicitly in terms of $F_\delta$ and $\O$.

First it is worth noting, with $\D'$ as in equation~(\ref{E8}),
\begin{eqnarray}
[\D '(a)]^\dag &=& [ (S \ot S) \D^T(S^{-1}(a))]^\dag 
=(\S^{-1} \ot \S^{-1}) \cdot T \cdot[\D(S^{-1}(a))^\dag] \no \\
&=&(\S^{-1} \ot \S^{-1}) \cdot T \cdot[\O \D(S^{-1}(a)^\dag) \O^{-1}] \no \\
&=&(\S^{-1} \ot \S^{-1}) [\O^T \D^T(\S(a^\dag)) (\O^T)^{-1}]. \no
\end{eqnarray}
Now using (\ref{E32})(ii) and (\ref{E33})(iv) respectively, we may write
\begin{eqnarray*}
\S(a)=w^{-1}S(a)w, \ \S^{-1}(a)=S^{-1}(w^{-1})S^{-1}(a)S^{-1}(w)
\end{eqnarray*}
so that
\begin{align}
[\D '(a)]^\dag&= W^{-1}
(S^{-1} \ot S^{-1}) [\O^T \D^T(w^{-1}S(a^\dag)w) (\O^T)^{-1}] W \no \\
&= W^{-1} (S^{-1} \ot S^{-1}) (\O^T)^{-1} 
\D_0 (S^{-1}(w) a^\dag S^{-1}(w^{-1}))(S^{-1} \ot S^{-1})\O^T 
W \no
\end{align}
where we have introduced the following operators
\begin{eqnarray*}
W = S^{-1}(w) \ot S^{-1}(w), \ W^{-1} = S^{-1}(w^{-1}) \ot S^{-1}(w^{-1})
\end{eqnarray*}
in order to simplify the notation. Thus, with $F_0$ as in equation~(\ref{E13}),
\begin{align}
[\D '(a)]^\dag&=
W^{-1} (S^{-1} \ot S^{-1}) (\O^T)^{-1} \cdot F_0 \D(S^{-1}(w))  \D(a^\dag) \D(S^{-1}(w^{-1}))F_0^{-1}
(S^{-1} \ot S^{-1})\O^T W. \label{S10}
\end{align}
On the other hand we have
\begin{eqnarray*}
[\D ' (a)]^\dag=[F_\delta \D(a) F_\delta^{-1}]^\dag=(F_\delta^{-1})^\dag \D(a)^\dag F_\delta^\dag = (F_\delta^{-1})^\dag \O \D(a^\dag) \O^{-1} F_\delta^\dag.
\end{eqnarray*}
By comparison with equation (\ref{S10}), it follows that the operator
\begin{equation}
\O^{-1} F_\delta^\dag W^{-1} (S^{-1} \ot S^{-1})(\O^T)^{-1} 
F_0 \D(S^{-1}(w)) \label{E36}
\end{equation}
must commute with the co-product $\D$. Below we show in fact that equation~(\ref{E36}) reduces
to $ 1 \ot 1$.

It is first useful to determine the behaviour of $\bar \g$ in equation~(\ref{E10})(ii)
under an arbitrary twist $G \in H \ot H$. Under the twisted structure induced by $G$ the
operator $\bar \g$ is twisted to $\bar \g_G$, given by equation~(\ref{E10})(ii,iii) for
the twisted structure, so that
\begin{align}
(i)&~~\bar \g_G= \bar A_i^G \b_G S(\bar D_i^G) \ot \bar B_i^G \b_G S( \bar C_i^G) \no \\
\textrm{where~~}(ii)&~~ \bar A_i^G \ot \bar B_i^G \ot \bar C_i^G \ot \bar D_i^G =
(\D_G\ot 1 \ot 1)\Phi_G^{-1}\cdot(\Phi_G \ot 1). \label{E25}
\end{align}
We have shown in a previous publication~\cite{mdgtl} that 
\begin{eqnarray}\label {P9}
\bar \g_G &=& G \ \D(g_i) \ \bar \g \ (S \ot S)(G^T\D^T(g^i)).
\end{eqnarray}

\begin{Proposition}\label{P16} Let $\gamma$ be the operator of equation~(\ref{E9})(ii). Then
\begin{eqnarray}
\gamma^\dag &=& \O \D(\O^i) \ (S^{-1} \ot S^{-1})\bar \gamma^T \ (S^{-1} \ot S^{-1})
[\O^T \D^T(\O_i)] \ W \no \\ 
&=&  \O \D(\O^i) \ (S^{-1} \ot S^{-1})\bar \gamma^T 
\ (S^{-1} \ot S^{-1}) \D^T(\O_i)
\ (S^{-1} \ot S^{-1})\O^T \ W. \no 
\end{eqnarray}
\end{Proposition}
\begin{proof}
       From equation~(\ref{E9})(ii) we have
       \begin{eqnarray}
       \gamma^\dag &=&  C_i^\dag \a^\dag \S^{-1}(B_i^\dag) 
       \ot D_i^\dag \a^\dag \S^{-1}(A_i^\dag) \no \\
       &\stackrel{(\ref{E32}')(ii)}{=}&  C_i^\dag \S^{-1}(w) \S^{-1}(\b_\O) \S^{-1}(B_i^\dag) 
       \ot D_i^\dag \S^{-1}(w) \S^{-1}(\b_\O) \S^{-1}(A_i^\dag) \no \\
       &\stackrel{(\ref{E33})(iii)}{=}&  C_i^\dag S^{-1}(w) \S^{-1}(\b_\O) \S^{-1}(B_i^\dag) 
       \ot D_i^\dag S^{-1}(w) \S^{-1}(\b_\O) \S^{-1}(A_i^\dag) \no \\
       &\stackrel{(\ref{E33})(iv)}{=}&
       [ C_i^\dag S^{-1}(\b_\O) S^{-1}(B_i^\dag) 
       \ot D_i^\dag  S^{-1}(\b_\O) S^{-1}(A_i^\dag)]\ W. \no
       \end{eqnarray}
       Now from equation~(\ref{E9})(iii)
       \begin{eqnarray}
        A_i^\dag \ot B_i^\dag \ot C_i^\dag \ot D_i^\dag &=&
       [(\Phi^{-1} \ot 1)(\D \ot 1 \ot 1) \Phi]^\dag \no \\
       &=& (\D_\O \ot 1 \ot 1) \Phi^\dag [(\Phi^\dag)^{-1} \ot 1] \no \\
       &=& (\D_\O \ot 1 \ot 1) \Phi_\O^{-1}(\Phi_\O \ot 1) \no \\
       &\stackrel{(\ref{E10})(iii)}{=}& 
        \bar{A}_i^\O \ot \bar{B}_i^\O \ot \bar{C}_i^\O \ot\bar{D}_i^\O \no  
       \end{eqnarray}
       which is the operator of equation~(\ref{E10})(iii) for the twisted structure
       induced by $\O$ (see equation~(\ref{E25})(ii)). Thus
       \begin{eqnarray}
       \gamma^\dag&=&[ \bar C_i^\O S^{-1}(\b_\O) S^{-1} (\bar B_i^\O) \ot 
       \bar D_i^\O S^{-1}(\b_\O) S^{-1}(\bar A_i^\O)] \ W \no \\
       &=&(S^{-1}\ot S^{-1})[ \bar B^\O_i \b_\O S(\bar C^\O_i) \ot \bar A^\O_i \b_\O
       S(\bar D^\O_i)]\ W \no \\
       &\stackrel{(\ref{E25})(i)}{=}&
       (S^{-1} \ot S^{-1})(\bar \gamma_\O^T) \ W \label{S11} \\
       &\stackrel{(\ref{P9})}{=}&(S^{-1} \ot S^{-1})\cdot T 
        \cdot [ \O \D(\O_i) \ \bar\gamma \ (S \ot S)(\O^T\D^T(\O^i))] \ W \no \\
       &=&\O \D(\O_i) \ (S^{-1} \ot S^{-1}) \bar\gamma^T \ (S^{-1} \ot S^{-1})(\O^T\D^T(\O^i)) \ W \no
       \end{eqnarray}
       where, as usual $\bar \gamma^T \equiv T \cdot \bar \gamma$. This proves the result.
\end{proof}

We are now in a position to compute $F_\delta^\dag$. From equation~(\ref{E9})(i) we have
immediately

\begin{eqnarray}
F_\delta^\dag & = & \D(\Z_\nu)^\dag \ \gamma^\dag \ (\S^{-1} \ot \S^{-1})\{\D^T(S^{-1} [\X_\nu \b S(\Y_\nu)])^\dag \} \no \\
              & = & \D(\Z_\nu)^\dag \ \gamma^\dag \ (\S^{-1} \ot \S^{-1}) \D_\O^T([S^{-1}(\X_\nu \b S(\Y_\nu))]^\dag )\no \\
              & = & \D_\O(\Z_\nu^\dag) \ \gamma^\dag \ (\S^{-1} \ot \S^{-1}) \D_\O^T(\S[\S^{-1}(\Y_\nu^\dag) \b^\dag \X_\nu^\dag])\no \\
	      & = & \D_\O(\Z_\nu^\dag) \ \gamma^\dag \ (\S^{-1} \ot \S^{-1}) \D_\O^T[\S(\X_\nu^\dag)\S(\b^\dag) \Y_\nu^\dag]\no \\
	      &\stackrel{(\ref{E32}')(i)}{=}& \D_\O(\Z_\nu^\dag) \ \gamma^\dag \ (\S^{-1} \ot \S^{-1}) \D_\O^T[\S(\X_\nu^\dag) w^{-1}\a_\O \Y_\nu^\dag]\no \\
	      &\stackrel{(\ref{E32})(ii)}{=}& \D_\O(\Z_\nu^\dag) \ \gamma^\dag \ (\S^{-1} \ot \S^{-1}) \D_\O^T[w^{-1} S(\X_\nu^\dag) \a_\O \Y_\nu^\dag]\no \\
	      & = & \D_\O(\Z_\nu^\dag) \ \gamma^\dag \ (\S^{-1} \ot \S^{-1}) \D_\O^T[S(\X_\nu^\dag) \a_\O \Y_\nu^\dag] \ (\S^{-1} \ot \S^{-1}) \D_\O^T(w^{-1})\no \\
	      &\stackrel{(\ref{S11})}{=}& \D_\O(\Z_\nu^\dag) \ (S^{-1} \ot S^{-1}) \bar\gamma^T_\O \ W 
	      \ (\S^{-1} \ot \S^{-1})\D_\O^T[S(\X_\nu^\dag \a_\O \Y_\nu^\dag)] \ (\S^{-1} \ot \S^{-1}) \D_\O^T(w^{-1})\no \\
	      &\stackrel{(\ref{E33})(iv)}{=}& \D_\O(\Z_\nu^\dag) \ (S^{-1} \ot S^{-1})\bar\gamma^T_\O \ (\S^{-1} \ot \S^{-1})\D_\O^T[S(\X_\nu^\dag) \a_\O \Y_\nu^\dag)] 
	      \ (S^{-1} \ot S^{-1}) \D_\O^T(w^{-1}) \ W \no \\
	      & = & (S^{-1} \ot S^{-1}) \cdot T \cdot [ \D_\O(S(\X_\nu^\dag) \a_\O \Y_\nu^\dag) \ \bar\gamma_\O \ (S \ot S )\D_\O^T(\Z_\nu^\dag)]
	      \ (S^{-1} \ot S^{-1}) \D_\O^T(w^{-1}) \ W. \no 
\end{eqnarray}
Now we may write
\begin{eqnarray*}
\X_\nu^\dag \ot \Y_\nu^\dag \ot \Z_\nu^\dag=(\Phi^\dag)^{-1}=\Phi_\O = X_\nu^\O \ot y_\nu^\O \ot Z_\nu^\O
\end{eqnarray*}
which is the co-associator for the twisted structure induced by $\O$. Thus we have,
\begin{eqnarray*}
F_\delta^\dag & = & (S^{-1} \ot S^{-1}) \cdot T \cdot [ \D_\O(S(X_\nu^\O) \a_\O Y_\nu^\O] \ \bar\gamma_\O \ (S \ot S )\D_\O^T(Z_\nu^\O)]
\ (S^{-1} \ot S^{-1}) \D_\O^T(w^{-1}) \ W \no \\
&\stackrel{(\ref{E10})(ii)}{=}&(S^{-1} \ot S^{-1}) \cdot T \cdot [(F_\delta^\O)^{-1}] \ (S^{-1} \ot S^{-1})\D_\O^T(w^{-1}) \ W \no 
\end{eqnarray*}
where $F_\delta^\O$ is the Drinfeld twist for the twisted structure induced by $\O$ and $(F_\delta^\O)^{-1}$ is its inverse. 

We now make use of the following theorem proved in~\cite{mdgtl}.
\begin{Theorem}\label{T4}
Let $G \in H \ot H$ be a twist on a QHA $H$. Then under the twisted structure induced by $G$, $F_\d^{-1}$ is twisted to
\begin{eqnarray*}
(F^G_\d)^{-1} \equiv (F^{-1}_\d)_G=G\cdot F_\d^{-1}\cdot(S \ot S)G^T.
\end{eqnarray*}
\end{Theorem}

It follows from Theorem \ref{T4} that
\begin{eqnarray}
F_\delta^\dag & = & (S^{-1} \ot S^{-1}) \cdot T \cdot [\O  F_\delta^{-1}  ( S \ot S) \O^T] \ (S^{-1} \ot S^{-1}) \D_\O^T(w^{-1}) \ W \no \\ 
              & = & \O (S^{-1} \ot S^{-1})(F_\delta^T)^{-1} \ ( S^{-1} \ot S^{-1})\O^T \ (S^{-1} \ot S^{-1})\D_\O^T(w^{-1}) \ W \no \\
              &\stackrel{(\ref{E13})}{=}& \O F_0^{-1} \ (S^{-1} \ot S^{-1})\D^T(w^{-1}) \ ( S^{-1} \ot S^{-1})\O^T \ W \no \\
              & = & \O  F_0^{-1}  \D_0(S^{-1}(w^{-1})) \ ( S^{-1} \ot S^{-1})\O^T \ W \no
\end{eqnarray}
with $\D_0$ the co-product of equation (\ref{E8}$'$). We thus arrive at our main result
\begin{Theorem}\label{T7}
\begin{eqnarray*}
F_\delta^\dag & = & \O \D(S^{-1}(w^{-1})) F_0^{-1} (S^{-1} \ot S^{-1})\O^T \ W.
\end{eqnarray*}
\end{Theorem}
\begin{Corollary}\label{T7C1}
With $F_0$ as in equation (\ref{E13}),
\begin{eqnarray*}
F_0^\dag & = & \O \D(w^{-1})F_\delta^{-1} \ (S \ot S)\O^T \ (w \ot w). 
\end{eqnarray*}
\end{Corollary}
\begin{proof}
       From equation (\ref{E13}), $F_0=(S^{-1} \ot S^{-1})F_\delta^T$, which implies
       \begin{eqnarray*}
       F_0^\dag=(\S \ot \S)[(F_\delta^T)^\dag].
       \end{eqnarray*}
       Now from Theorem \ref{T7} we get
       \begin{eqnarray}
       F_0^\dag & = & (\S \ot \S)[\O^T \D^T(S^{-1}(w^{-1})) \ (F_0^T)^{-1}] \ (S^{-1} \ot S^{-1})\O \ W \no \\
                & = & (\S \ot \S) W (\S \ot \S)[\O^T \D^T(S^{-1}(w^{-1})) (F_0^T)^{-1} \ (S^{-1} \ot S^{-1})\O] \no \\
                &\stackrel{(\ref{E33})(iii)}{=}& (w \ot w) (\S \ot \S)[\O^T \D^T(S^{-1}(w^{-1})) (F_0^T)^{-1}  (S^{-1} \ot S^{-1})\O] \no \\
                &\stackrel{(\ref{E32})(ii)}{=}& (S \ot S)[\O^T \D^T(S^{-1}(w^{-1}))(F_0^T)^{-1} (S^{-1} \ot S^{-1})\O] (w \ot w). \no 
       \end{eqnarray}
       Now from equation (\ref{E13}), $(F_0^T)^{-1}=(S^{-1} \ot S^{-1})(F_\delta^{-1})$,
       so that
       \begin{eqnarray}
       F_0^\dag & = & (S \ot S) [\O^T \D^T(S^{-1}(w^{-1})) (S^{-1} \ot S^{-1})F_\delta^{-1}(S^{-1} \ot S^{-1})\O] (w \ot w) \no \\
                & = & \O F_\delta^{-1} \D'(w^{-1}) (S \ot S)\O^T (w \ot w) \no \\
                & = & \O \D(w^{-1}) F_\delta^{-1} (S \ot S)\O^T (w \ot w) \no
       \end{eqnarray}
       which proves the result.
\end{proof}
\begin{Corollary}\label{T7C2} The operator of equation (\ref{E36}) is given by
\begin{eqnarray*}
\O^{-1} F_\delta^\dag W (S^{-1} \ot S^{-1})(\O^T)^{-1}  F_0  \D(S^{-1}(w))= 1 \ot 1 .
\end{eqnarray*}

\end{Corollary}
\setcounter{Corollary}{0}
\begin{proof}
       Follows by an easy computation using Theorem \ref{T7}.
\end{proof}

If $S$ is $*$-compatible, so that $\S=S$, the above result for $F_\delta^\dag$
remains unaltered, except for the simplification that $w$ is central.

The results above have a number of interesting consequences. In particular, we are now in
a position to obtain the $*$-canonical elements of equation (\ref{E35}) pertinent to the
$*$-QHAs of propositions \ref{P2}, \ref{P2}$'$. By a straightforward calculation using 
Theorem \ref{T7} and Corollary \ref{T7C1}, we immediately obtain
\begin{Proposition}\label{P17}
\begin{eqnarray*}
\O'&\equiv&(F_\delta^\dag)^{-1}\O F_\delta^{-1}\\
& = & W (S^{-1} \ot S^{-1})(\O^T)^{-1} F_0 \D(S^{-1}(w))F_\delta^{-1}  \\
& = & [S^{-1}(w^{-1}) \ot S^{-1}(w^{-1})]  (S^{-1} \ot S^{-1})(\O^T)^{-1} F_0 \D(S^{-1}(w))F_\delta^{-1}  \\
\O_0&\equiv&(F_0^\dag)^{-1}\O F_0^{-1}\\
&=&(w^{-1} \ot w^{-1})  (S \ot S)(\O^T)^{-1} F_\delta \D(w) F_0^{-1}.  
\end{eqnarray*}
\end{Proposition}

All of the results of this section will obviously hold with $\O$ replaced
by $\O^\dag$ in which case $w$ must be replaced by $\bar w$. 
In particular the result of proposition \ref{P16} will hold with $\O$ replaced
by $\O^\dag$ and $w$ by $\bar w$. Taking the Hermitian conjugate of the resultant
expression, using proposition \ref{P15}, it is then easy to obtain an expression
for $\bar \gamma^\dag$ in terms of $\O,w$ and $\gamma$.

Replacing $\O,w$ with $\bar \O, \bar w$ respectively in proposition \ref{P17},
we arrive at the corresponding conjugate $*$-canonical elements
\vskip 3mm \noindent
{\bf Proposition \ref{P17}$'$~~}{\it
\begin{eqnarray*}
(\O')^\dag&=&(F_\delta^\dag)^{-1}\O^\dag F_\delta^{-1}\\
& = & [S^{-1}(\bar w^{-1}) \ot S^{-1}(\bar w^{-1})] (S^{-1} \ot S^{-1}){(\O^\dag)^T}^{-1} F_0 \D(S^{-1}(\bar w))F_\delta^{-1}  \\
\O_0^\dag&=&(F_0^\dag)^{-1}\O^\dag F_0^{-1}\\
& = & [\bar w^{-1} \ot \bar w^{-1}] (S \ot S){(\O^\dag)^T}^{-1} F_\delta \D(\bar w) F_0^{-1}.  
\end{eqnarray*}
}\\
{\it Note} Replacing $\O$ with $\O^\dag$ and $w$ with $\bar w$ in Corollary \ref{T7C1} to
Theorem \ref{T7} gives
\begin{eqnarray*}
F_0^\dag & = & \O \D(w^{-1}) F_\delta^{-1} (S \ot S)\O^T (w\ot w) \\
         & = & \O^\dag \D(\bar w^{-1}) F_\delta^{-1} (S \ot S)(\O^T)^\dag (\bar w \ot \bar w)
\end{eqnarray*}
which implies that
\begin{eqnarray*}
\O^{-1} \O^\dag \D(\bar w^{-1}) F_\delta^{-1} (S \ot S)(\O^T)^\dag  = \D(w^{-1}) F_\delta^{-1} (S \ot S)\O^T (c^{-1} \ot c^{-1})
\end{eqnarray*}
where $c=w^{-1} \bar w$ is the central element of equation (\ref{E34}). Using the fact that
$\O^{-1}\O^\dag$ commutes with $\D$, then gives
$$(\O^{-1} \O^\dag) F_\delta^{-1} (S \ot S)(\O^T)^\dag (S \ot S)(\O^T)^{-1} = \D(c)F_\delta^{-1} (c^{-1} \ot c^{-1})$$
or, to put it another way
$$(\O^{-1} \O^\dag) F_\delta^{-1} (S \ot S)(\O^{-1} \O^\dag)^T F_\delta= (c^{-1} \ot c^{-1}) \D(c).$$
In the case $H$ is self conjugate, so that $\O=\O^\dag$, this just reduces to an identity.

We now consider the important case when the $*$-QHA $H$ is quasi-triangular.


\section{The Quasi-triangular Case}\label{starQTQHA}

A quasi-Hopf algebra $H$ is called quasi-triangular if there exists an invertible element $\R \in H\ot H$
called the $R$-matrix, such that
\begin{align}
	\D^T(a)\R & = \R\D(a),\quad \forall a\in H \label{intertwine} \\
	{(\D\ot 1)}\R & = \Phi_{231}^{-1} \R_{13}\Phi_{132} \R_{23}\Phi_{123}^{-1} \label{dot1} \\
	{(1\ot\D)}\R & = \Phi_{312} \R_{13}\Phi_{213}^{-1} \R_{12}\Phi_{123}.\label{1otd}
\end{align}

Above, $\D^T(a)=T \cdot \D(a)$, where $T:H \ot H \rightarrow H \ot H$ is the usual twist map, $T(a \ot b) = b \ot a$. 
For the co-associator we have followed the conventions of~\cite{cas,oqhsa} so that,
\begin{equation*}
	\Phi_{231} = \sum_\nu Z_\nu \ot Y_\nu \ot X_\nu,
	\quad \Phi^{-1}_{312} = \sum_\nu \bar Z_\nu \ot 
	\bar X_\nu \ot \bar Y_\nu,
	\quad \text{etcetera.}
\end{equation*}
We set 
\begin{equation*}
	 \R = \sum_i e_i \ot e^i, \ \R^{-1} = \sum_i \bar{e}_i \ot \bar{e}^i
\end{equation*}
in terms of which  
\begin{equation*}
	\R_{12}= \sum_i e_i \ot e^i \ot 1,\quad \R_{13} = \sum_i e_i \ot 1 \ot e^i, \quad \text{etcetera.}
\end{equation*}

Throughout this section we assume $H$ is a $*$-QHA with $*$-canonical element $\O$, which is moreover quasi-triangular,
i.e. admits an $R$-matrix, $\R \in H \ot H$ satisfying equations (\ref{intertwine}-\ref{1otd}).

The $R$-matrix $\R$ satisfies the additional relations 
\begin{equation*}
	(\epsilon \otimes 1) \R = (1 \otimes \epsilon) \R = 1,
\end{equation*}
which follow from (\ref{dot1}) and (\ref{1otd}). Since $\R$ is invertible and satisfies the co-unit property~(\ref{twist-counit}) it qualifies as a twist.

Recall that if $H$ is a QTQHA then it is also a QTQHA under the opposite structure of proposition~\ref{P1} but with opposite $R$-matrix $\R^T \equiv T \cdot \R$. 
Twisting $H$ with the $R$-matrix $\R$ gives rise to this opposite structure but with antipode $S$ and canonical elements $\a_{\R}, \b_{\R}$ given by equ.~(\ref{TWCAN}). Applying theorem~\ref{T1} to these equivalent QHA structures gives 

\begin{Theorem}
There exists a unique invertible $u \in H$ such that
\begin{equation*}
S(a)=uS^{-1}(a)u^{-1}, \textrm{~or~} S^2(a)=uau^{-1},\quad \forall a \in H
\end{equation*}
and
\begin{equation}
uS^{-1}(\a)=\a_\R,\quad\b_\R u=S^{-1}(\b). \label{E18}
\end{equation}
Explicitly,
\begin{eqnarray}
u&=& S(Y_\nu\b S(Z_\nu))\a_\R X_\nu =
 S(\Z_\nu) \a_\R \Y_\nu S^{-1}(\b)S^{-1}(\X_\nu)\no \\
u^{-1}&=&  Z_\nu \b_\R S(S(X_\nu)\a Y_\nu) =
 S^{-1}(\Z_\nu)S^{-1}(\a) \Y_\nu \b_\R S(\X_\nu). \label{E19}
\end{eqnarray}
\end{Theorem}

Now from the intertwining property (\ref{intertwine}) we have
\begin{equation*}
\R \D(a^\dag)=\D^T(a^\dag) \R
\end{equation*}
which gives, upon applying $\dag$,
$$\R^\dag \tilde{\D}^T(a)=\tilde{\D}(a)\R^\dag,
~\textrm{~or\quad} (\R^\dag)^{-1}\tilde{\D}(a)=\tilde{\D}^T(a)(\R^\dag)^{-1},\quad\forall a \in H. $$ 
Thus $(\R^\dag)^{-1}$ satisfies the intertwining property (\ref{intertwine}) for the co-product
$\tilde{\D}$ of equation (\ref{E28}), given by [cf equation(\ref{E29}$'$)]
$$\tilde{\D}(a)=\D(a^\dag)^\dag=\D_\O(a)=\D_{\O^\dag}(a),\quad\forall a \in H.$$
We thus expect $(\R^\dag)^{-1}$ to give rise to an $R$-matrix for $H$ with the QHA
structure of proposition \ref{P11}, which is indeed the case.
\vskip 3mm \noindent
{\bf Proposition \ref{P11}$'$~~}{\it
Suppose $H$ is any quasi-triangular QHA admitting a conjugation operation $\dag$ satisfying
only equation (\ref{EADAG}). Then $H$ is also a quasi-triangular QHA with the structure
of proposition \ref{P11} with $R$-matrix $(\R^\dag)^{-1}$.
}
\noindent
\begin{proof}
       First recall that, with the structure of proposition \ref{P11}, $H$ is a QHA
       with the same co-unit but with co-product $\tilde{\D}$, co-associator $\tilde \Phi =
       (\Phi^\dag)^{-1}$ and antipode $\S$ given by equation (\ref{E30}). We have 
       already seen that if $\R$ is an $R$-matrix for $H$ then $(\R^\dag)^{-1}$
       satisfies the intertwining property (\ref{intertwine}) for this structure. It thus
       remains to consider (\ref{dot1}, \ref{1otd}).
       
       Taking the conjugate inverse of equations (\ref{dot1}, \ref{1otd}) and using
       $\dag\cdot\D=\tilde\D\cdot\dag$, gives immediately
       \begin{eqnarray*}
       (\tilde \D \ot 1)(\R^\dag)^{-1} & = & \Phi^\dag_{231}(\R^\dag)^{-1}_{13}(\Phi^\dag)^{-1}_{132}(\R^\dag)^{-1}_{23}\Phi^\dag_{123} \\
       (1 \ot \tilde \D )(\R^\dag)^{-1} & = & 
       (\Phi^\dag)^{-1}_{312}(\R^\dag)^{-1}_{13}\Phi^\dag_{213}(\R^\dag)^{-1}_{12}(\Phi^\dag)^{-1}_{123}.
       \end{eqnarray*}
       Setting $\tilde \Phi = (\Phi^\dag)^{-1}$, which is the co-associator for this
       structure, implies
       \begin{eqnarray*}
       (\tilde \D \ot 1)(\R^\dag)^{-1} & = & \tilde\Phi^{-1}_{231}(\R^\dag)^{-1}_{13} \tilde 
       \Phi_{132}(\R^\dag)^{-1}_{23}\tilde\Phi^{-1}_{123} \\
       (1 \ot \tilde \D )(\R^\dag)^{-1} & = & \tilde\Phi_{312}(\R^\dag)^{-1}_{13} \tilde \Phi^{-}_{213}(\R^\dag)^{-1}_{12}\tilde\Phi_{123}
       \end{eqnarray*}
       as required.
\end{proof}
\begin{Corollary*}
In the case $H$ is a $*$-QHA with $*$-canonical element $\O$,
\begin{equation}
\bar \R = (\O^T)^{-1}(\R^\dag)^{-1}\O \label{E37}
\end{equation}
determines an $R$-matrix for $H$
\end{Corollary*}
\begin{proof}
       In such a case the QBA structure of proposition \ref{P11} is obtained by
       twisting with $\O$ (or $\O^\dag$); i.e. $\tilde{\D}=\D_\O, \tilde \Phi = \Phi_\O$.
       The result above shows that $(\R^\dag)^{-1}$ is an $R$-matrix for this twisted
       structure. Since $\bar \R$ of equation (\ref{E37}) is obtained from $(\R^\dag)^{-1}$
       by twisting with $\O^{-1}$ it follows that $\bar \R$ must determine an $R$-matrix
       for $H$, i.e. satisfy equations (\ref{intertwine} - \ref{1otd}), since $\O^{-1}$
       will ``undo'' the twist $\O$.
\end{proof}
\begin{remark*}
The above corollary may be proved directly.
\end{remark*}

Thus associated with the $R$-matrix $\bar \R$ of equation (\ref{E37}) we have a $u$-operator
$\bar u = u_{\bar \R}$ and its inverse $\bar u^{-1}$ given explicitly by equation
(\ref{E19}) with $\R$ replaced by $\bar \R$. Then $\bar u \in H$ is the unique operator
satisfying
\begin{equation}
S^2(a)=\bar u a \bar u^{-1},~\forall a \in H;\quad\bar u S^{-1}(\a)=\a_{\bar \R},
\quad \b_{\bar R} \bar u = S^{-1}(\b). \label{E38}
\end{equation}
Here we explore the connection between $u$ and $\bar u$. We first need, with $\S$ as in
equation (\ref{E30}), the following
\begin{Lemma}\label{L10}(notation as in lemma \ref{L8})
\begin{eqnarray*}
&(i)&\S^{-1}(\O_i)\b_\R^\dag\O^i=S^{-1}(w^{-1})S^{-1}(\a_{\bar \R})  \\
&(ii)&\bar \O_i\a_\R^\dag\S^{-1}(\bar \O^i)=S^{-1}(\b_{\bar \R})S^{-1}(w). 
\end{eqnarray*}
\end{Lemma}
\begin{proof}
       (i) First note, from equation (\ref{E37}), that $(\R^\dag)^{-1}=\O^T \bar \R \O^{-1}$. Hence
       using lemma \ref{L9}(iii) we have
       \begin{eqnarray*}
       \S^{-1}(\O_i)\b_\R^\dag\O^i&=&\S^{-1}(\O_i)S^{-1}(w^{-1})S^{-1}[\a_{(\R^\dag)^{-1}\O}]\O^i  \\
       &=&\S^{-1}(\O_i)S^{-1}(w^{-1})S^{-1}[\a_{\O^T \bar \R}]\O^i  \\
       &\stackrel{(\ref{E33})(iv)}{=}& S^{-1}(w^{-1})S^{-1}(\O_i)S^{-1}[\a_{\O^T \bar \R}]\O^i  \\
       &=&S^{-1}(w^{-1})S^{-1}[ S(\O^i)\a_{\O^T \bar \R}\O_i]  \\
       &=&S^{-1}(w^{-1})S^{-1}(\a_{\bar \R}). 
       \end{eqnarray*}
       Part (ii) is proved in a similar fashion. 
\end{proof}
We are now in a position to compute $u^\dag$. From equation (\ref{E19}) we have immediately
\begin{equation}
u^\dag= \S(\X_\nu^\dag) S^{-1}(\b)^\dag \Y_\nu^\dag \a_\R^\dag \S^{-1}(\Z_\nu^\dag), \label{S12}
\end{equation}
where
\begin{eqnarray}
\X_\nu^\dag \ot \Y_\nu^\dag \ot \Z_\nu^\dag &=& (\Phi^\dag)^{-1} \no \\
&=&\Phi_\O = (\O \ot 1)\cdot(\D \ot 1)\O\cdot\Phi\cdot(1\ot\D)\O^{-1}\cdot(1\ot\O^{-1}) \no \\ 
&=&\O_i\O_j^{(1)}X_\nu\bar\O_k \ot \O^i\O_j^{(2)}Y_\nu\bar\O^k_{(1)}\bar\O_l \ot \O^jZ_\nu\bar\O^k_{(2)}\bar\O^l \no
\end{eqnarray}
where we have adopted the obvious notation (all repeated indices to be summed over). Substituting
into~(\ref{S12}) then gives
\begin{eqnarray*}
u^\dag&=&  \S(\O_i\O_j^{(1)}X_\nu\bar\O_k) \S^{-1}(\b)^\dag 
\ \O^i\O_j^{(2)}Y_\nu\bar\O^k_{(1)}\bar\O_l \a_\R^\dag \S^{-1}(\O^jZ_\nu \bar\O^k_{(2)}\bar\O^l)  \\
&\stackrel{(\ref{E32})(i)}{=}& 
\S(\O_i\O_j^{(1)}X_\nu\bar\O_k) w^{-1} \a_\O \O^i\O_j^{(2)}Y_\nu\bar\O^k_{(1)}\bar\O_l \a_\R^\dag \S^{-1}(\bar\O^l)\S^{-1}(\O^jZ_\nu\bar\O^k_{(2)})  \\
&\stackrel{(\ref{E32})(ii)}{=}& w^{-1}
S(\O_i\O_j^{(1)}X_\nu\bar\O_k) \a_\O \O^i\O_j^{(2)}Y_\nu\bar\O^k_{(1)}\bar\O_l \a_\R^\dag \S^{-1}(\bar\O^l)\S^{-1}(\O^jZ_\nu\bar\O^k_{(2)}) \no \\
&=& w^{-1}
S(\O_j^{(1)}X_\nu\bar\O_k)S(\O_i) \a_\O \O^i\O_j^{(2)}Y_\nu\bar\O^k_{(1)}\bar\O_l \a_\R^\dag \S^{-1}(\bar\O^l)\S^{-1}(\O^jZ_\nu\bar\O^k_{(2)}). \no
\end{eqnarray*} 
Now,
\begin{equation*}
	S(\O_i) \a_\O \O^i=(\a_\O)_{\O^{-1}}=\a_{\O^{-1}\O}=\a
\end{equation*}
and from lemma \ref{L10}(ii)
\begin{eqnarray*}
\bar\O_l \a_\R^\dag \S^{-1}(\bar\O^l)&=& S^{-1}(\b_{\bar \R})S^{-1}(w)\stackrel{(\ref{E38})}{=}
S^{-1}[S^{-1}(\b)\bar u^{-1}]\cdot S^{-1}(w)  \\
&=&S^{-1}[\bar u^{-1} S(\b)]S^{-1}(w)=\b S^{-1}(\bar u^{-1})S^{-1}(w). 
\end{eqnarray*}
At this point it is worth noting that,
\begin{equation}
S^{-1}(\bar u^{-1})S^{-1}(a)S^{-1}(\bar u)=S^{-1}[\bar u a \bar u^{-1}]\stackrel{(\ref{E38})}{=}
S^{-1}[S^2(a)]=S(a).\label{S13}
\end{equation}
Substituting into the above gives
\begin{eqnarray*}
u^\dag&=& w^{-1} S(\O_j^{(1)}X_\nu\bar\O_k) \a \O_j^{(2)} Y_\nu \bar\O^k_{(1)} \b 
\ S^{-1}(\bar u^{-1})S^{-1}(w)\S^{-1}(\O^jZ_\nu\bar\O^k_{(2)}) \no \\
&\stackrel{(\ref{E33})(iv)}{=}& w^{-1}
S(\O_j^{(1)}X_\nu\bar\O_k) \a \O_j^{(2)} Y_\nu \bar\O^k_{(1)} \b 
 \ S^{-1}(\bar u^{-1})S^{-1}(\O^jZ_\nu\bar\O^k_{(2)})S^{-1}(w) \no \\
&\stackrel{(\ref{S13})}{=}& w^{-1}
S(\O_j^{(1)}X_\nu\bar\O_k) \a \O_j^{(2)} Y_\nu \bar\O^k_{(1)} \b \ S(\O^jZ_\nu\bar\O^k_{(2)})S(w) S^{-1}(\bar u^{-1})\no \\
&=& w^{-1}
S(X_\nu\bar\O_k) S(\O_j^{(1)}) \a \O_j^{(2)} Y_\nu \bar\O^k_{(1)} \b \ S(\bar\O^k_{(2)}) S(\O^jZ_\nu)S(w) S^{-1}(\bar u^{-1})\no \\
&\stackrel{(\ref{Sab})}{=}& w^{-1}
S(X_\nu) \a Y_\nu \b S(Z_\nu)S(w) S^{-1}(\bar u^{-1})\no \\
&\stackrel{(\ref{Sab})}{=}& 
w^{-1} S(w) S^{-1}(\bar u^{-1}). \no
\end{eqnarray*}
Hence we have proved
\begin{Proposition}\label{P18}
$$u^\dag=w^{-1} S(w) S^{-1}(\bar u^{-1}).$$
\end{Proposition}

With regard to proposition~\ref{P18} it is worth noting the following result concerning the
antipode $\S$ of equation~(\ref{E30}):
\begin{Lemma}\label{L11}
$v = w^{-1} S(w) u$ determines a $u$-operator for the quasi-triangular QHA structure of proposition~\ref{P11}$\ '$, i.e.
$$\S(a)=v\S^{-1}(a)v^{-1},\quad\forall a \in H.$$
\end{Lemma}
\begin{proof}
       From equation~(\ref{E32})(ii) we have
       $$w\S(a)w^{-1}=S(a)=uS^{-1}(a)u^{-1}
       \stackrel{(\ref{E33})(iv)}{=}uS^{-1}(w)\S^{-1}(a)S^{-1}(w^{-1})u^{-1}$$
       which implies
       $$\S(a)=w^{-1}u S^{-1}(w)\S^{-1}(a)S^{-1}(w^{-1})u^{-1}w,\quad\forall a \in H$$
       so that
       $$v=w^{-1}uS^{-1}(w)=w^{-1}S(w)u$$
       is a $u$-operator for $\S$.
\end{proof}
\begin{Corollary*}
$u^\dag$ is also a $u$-operator for $\S$.
\end{Corollary*}
\begin{proof}
       This follows from lemma~\ref{L11} and proposition~\ref{P18} by noting from equation~(\ref{S13}) above,
       that $S^{-1}(\bar u^{-1})$ is a $u$-operator for $H$ (with respect to $S$).
\end{proof}

Following the definition (\ref{EADAG}-\ref{PHIDAG}) of a $*$-QHA, it is natural to define a quasi-triangular $*$-QHA ($*$-QTQHA)
as one for which the complete QHA structure of proposition~\ref{P11}$'$ is obtainable [modulo $(S,\a,\b)$] by twisting with $\O$.
However, since $(\R^T)^{-1}$ is also an $R$-matrix this leads to two natural classes of $*$-QTQHA:
\begin{Definition}
A $*$-QHA with $*$-canonical element $\O$ is called a $*$-QTQHA of Type I (resp. Type II) if it is a
quasi-triangular QHA with $R$-matrix satisfying
\begin{equation}
(\R^\dag)^{-1}= \O^T \R \O^{-1},\quad [\textrm{resp.~} \O^T (\R^T)^{-1} \O^{-1}]. \label{E39}
\end{equation}
\end{Definition}

In the Hopf algebra setting Majid~\cite{majid} has pointed out that two natural Hopf-$*$
structures arise in the quasi-triangular case, called the antireal case, $\R^\dag=\R^{-1}$, 
and the real case $\R^\dag=\R^T$. For $*$-QTQHA we see that upon setting $\O= 1 \ot 1$
the type I case reduces to Majid's antireal case and the type II to his real case.
In the triangular case, corresponding to $\R^T\R=1\ot 1$, or $\R^T=\R^{-1}$, the type I and 
type II cases coincide.

Our main result here is,
\begin{Theorem}[Twist Invariance]\label{T8}
Let $F \in H\ot H$ be an arbitrary twist on a $*$-QTQHA of type I (resp. type II). Then $H$ is also a
$*$-QTQHA of type I (resp. type II) under the twisted structure induced by $F$.
\end{Theorem}
\begin{proof}
       Following Theorem \ref{T5}, it suffices to prove that under the twisted structure induced
       by $F$ equation (\ref{E39}) holds: recall that $H$ is also a quasi-triangular QHA under this 
       structure with $R$-matrix $\R_F=F^T \R F^{-1}$. For the type I case we have
       \begin{eqnarray*}
       (\R_F^\dag)^{-1}&=&(F^T \R F^{-1})^{\dag^{-1}}  \\
       &=&(F^T)^{\dag^{-1}}(\R^\dag)^{-1}F^\dag \\
       &=&(F^T)^{\dag^{-1}}\O^T\R\O^{-1}F^\dag \\
       &=&(F^T)^{\dag^{-1}}\O (F^T)^{-1} \R_F F\O^{-1}F^\dag \\
       &=&\O^T_F \R_F \O^{-1}_F 
       \end{eqnarray*}
       where $\O_F=(F^\dag)^{-1}\O F^{-1}$ is the $*$-canonical element for the twisted structure,
       which shows that $H$ is also a type I $*$-QTQHA under this structure. The proof is similar for the 
       type II case.
\end{proof}

Thus in the type I case the $R$-matrix $\bar \R$ of equation (\ref{E37}) is given by
$$\bar \R =(\O^T)^{-1}(\R^\dag)^{-1}\O=(\O^T)^{-1}(\O^T R\O^{-1})\O=R$$
while in the type II case
$$\bar \R =(\O^T)^{-1}(\R^\dag)^{-1}\O=(\O^T)^{-1}(\O^T(\R^T)^{-1}\O^{-1})\O=(\R^T)^{-1}.$$
Thus for the $u$-operator $\bar u$ of equation (\ref{E38}) we have
$$\bar u=\left\{ 
\begin{array}{ll}
u  & \textrm{~~~~~(type I case)} \\
\tilde u & \textrm{~~~~~(type II case)}
\end{array}
\right.$$
where 
$$\tilde u = S(Y_{\nu}\b S(Z_\nu)) \a_{\tilde \R} X_\nu $$ which was shown in~\cite{mdgtl} to be 
given by $\tilde u =S(u^{-1})$. In view of proposition~\ref{P18} we thus arrive at 
\vskip 3mm \noindent
{\bf Proposition \ref{P18}$'$} {\it Let $H$ be a $*$-QTQHA. Then the $u$-operator of equation (\ref{E19})
must satisfy
$$u^\dag=w^{-1}S(w)\cdot \left\{ 
\begin{array}{ll}
S(u^{-1})  & \textrm{~~~~~(type I case)} \\
u & \textrm{~~~~~(type II case)}.
\end{array}
\right.$$ }
\noindent
For the type I case above, we used the well known result $S(u^{-1})=S^{-1}(u^{-1})$, as is easily
verified.

It is easily verified that 
\begin{equation}
z_u \equiv uS(u)=S(u)u \label{E40}
\end{equation}
is a central element, as shown in \cite{cas}. In terms of $z_u$ the result of proposition (\ref{P18}$'$) is
expressible as 
$$u^\dag=w^{-1}S(w)\cdot \left\{ 
\begin{array}{ll}
z_u^{-1}u  & \textrm{~~~~~(type I case)} \\
z_u S(u^{-1}) & \textrm{~~~~~(type II case)}.
\end{array}
\right.$$
The main difference between the type I and type II $*$-cases lies in the nature of the central
element $z_u$: it is always unitary in the type I case while in the type II case it is self-adjoint. 
Explicitly
\begin{Lemma}\label{L12}
Let $H$ be a $*$-QTQHA and $z_u$ the central element of equation (\ref{E40}). Then
$$z_u^\dag=\left\{ 
\begin{array}{ll}
z_u^{-1}  & \textrm{~~~~~(type I case)} \\
z_u  & \textrm{~~~~~(type II case)}.
\end{array}
\right.$$
\end{Lemma}
\begin{proof}
       First observe from equation (\ref{E33})(iv) that
       \begin{align}
       \S^{-1}(a)=S^{-1}(w^{-1})S^{-1}(a)S^{-1}(w),\quad\forall a \in H. \tag{$*$}
       \end{align}
       Now, with $z_u=uS(u)$ we have
       $$z_u^\dag=\S^{-1}(u^\dag)u^\dag\stackrel{(*)}{=}S^{-1}(w^{-1})S^{-1}(u^\dag)S^{-1}(w)u^\dag. $$
       Thus in the type I case 
       \begin{eqnarray}
       z_u^\dag&=&S^{-1}(w^{-1})S^{-1}[w^{-1}S(w)S(u^{-1})]\cdot S^{-1}(w)w^{-1}S(w)S(u^{-1}) \no \\
       &=&S^{-1}(w^{-1})[u^{-1}wS^{-1}(w^{-1})]\cdot S^{-1}(w)w^{-1}S(w)S(u^{-1}) \no \\
       &=&S^{-1}(w^{-1})u^{-1}S(w)S(u^{-1}) \no \\
       &=&S^{-1}(w^{-1})S^{-1}(w)u^{-1}S(u^{-1}) \no \\
       &=&u^{-1}S(u^{-1})=z_u^{-1},\no
       \end{eqnarray}
       while in the type II case 
       \begin{eqnarray}
       z_u^\dag&=&S^{-1}(w^{-1})S^{-1}[w^{-1}S(w)u]\cdot S^{-1}(w)w^{-1}S(w)u \no \\
       &=&S^{-1}(w^{-1})[S^{-1}(u)wS^{-1}(w^{-1})]\cdot S^{-1}(w)w^{-1}S(w)u \no \\
       &=&S^{-1}(w^{-1})S^{-1}(u)S(w)u \no \\
       &=&S^{-1}(w^{-1})S^{-1}(w)S^{-1}(u)u \no \\
       &=&S^{-1}(u)u=uS(u)=z_u. \no
       \end{eqnarray}  
\end{proof}

Lemma~\ref{L12} holds quite generally, regardless of whether or not $\O$ is self-adjoint or $S$ is $*$-compatible. It is a
universal property completely independent of $\O$.

Theorem~\ref{T8} shows that the category of type I or type II $*$-QTQHAs is invariant under twisting.
Now $H$ is also a quasi-triangular QHA with $R$-matrix $\R^T$ under the opposite structure of proposition~\ref{P1} and is
obtainable by twisting with $\R$. Moreover, proposition~\ref{P12} shows that $H$ is also a $*$-QHA
under this opposite structure with $*$-canonical element $\O^T = T \cdot \O$. It is therefore not 
surprising that we have the following extension
\vskip 3mm \noindent
{\bf Proposition \ref{P12}$'$}{\it~ A type I (resp. type II) $*$-QTQHA is also a type I (resp. type II)
$*$-QTQHA under the opposite structure of proposition~\ref{P12} with $R$-matrix $\R^T=T\cdot\R$}.
\noindent
\begin{proof}
       In view of the above and Proposition~\ref{P12} it remains to check equation (\ref{E39}) for
       this opposite structure. To this end we have in the type I case
       $$(\R^T)^{\dag^{-1}}\stackrel{(\ref{E39})}{=}\O\R^T(\O^T)^{-1}=(\O^T)^T\R^T(\O^T)^{-1}$$
       and similarly for the type II case
       $$(\R^T)^{\dag^{-1}}\stackrel{(\ref{E39})}{=}\O\R^{-1}(\O^T)^{-1}=(\O^T)^T{(\R^T)^T}^{-1}(\O^T)^{-1}$$
       which proves equation (\ref{E39}) for the opposite structure as required.
\end{proof}
It is important to note that the definition of $*$-QTQHA depends explicitly on the $*$-canonical
element $\O$ which is interconnected with the $R$-matrix through equation (\ref{E39}). Indeed,
if $\O_1=\O C$ is another $*$-canonical element with $C$ a compatible twist 
[see Theorem~\ref{T6}] then $H$ will not generally be a $*$-QTQHA with respect to $\O_1$ as is
easily seen. However, following Theorem~\ref{T8}, $H$ will be a $*$-QTQHA with twisted canonical
element $\O_C=(C^\dag)^{-1} \O C^{-1}$ and $R$-matrix $\R_C=C^T\R C^{-1}$.

As noted above, the definition of a $*$-QTQHA depends on the $*$-canonical element $\O$ (as well as $\R$). We
in fact have the following extension of Theorem~\ref{T6}:
\vskip 3mm \noindent
{\bf Theorem \ref{T6}$'$}{\it~ Let $H$ be a $*$-QTQHA with $*$-canonical element $\O$ and $R$-matrix $\R$. Then
$H$ is also a $*$-QTQHA with the same $R$-matrix $\R$ but with $*$-canonical element $\Gamma$ if and only
if there exists a compatible twist $C \in H\ot H$ such that $\Gamma=\O C$ and}
\begin{equation}
C^T \R C^{-1}=\R. \label{E41}
\end{equation}
\begin{proof}
       First from Theorem~\ref{T6}, in order for $\Gamma$ to be a $*$-canonical element there must
       exist a compatible twist $C \in H\ot H$ such that $\Gamma=\O C$. Now suppose $H$ is a $*$-QTQHA of
       type I, so that 
       $$(\R^\dag)^{-1}=\O^T \R \O^{-1}. $$
       Then in order for $H$ to be a $*$-QTQHA of type I with respect to $\Gamma$ it is necessary and sufficient
       that
       \begin{align}
       (\R^\dag)^{-1}=\Gamma^T\R\Gamma^{-1}  
       &\Leftrightarrow \Gamma^T\R \Gamma^{-1}=\O^T \R \O^{-1}  \no \\ 
       &\Leftrightarrow  \O^TC^T\R C^{-1}\O^{-1}=\O^T\R\O^{-1} \no \\
       &\Leftrightarrow C^T \R C^{-1}=\R
       \end{align}
       and similarly for the type II case. This proves the result.
\end{proof}
We thus have
\begin{Definition} We call a compatible twist $C$ on a quasi-triangular QHA
a quasi-triangul\-ar compatible twist if it satisfies equation (\ref{E41}).
\end{Definition}

Twisting a quasi-triangular QHA with such a twist will leave the entire structure unchanged
(modulo $\a,\b$). 
Quasi-triangular compatible twists on a quasi-triangular QHA $H$ form a subgroup of the group
of compatible twists on $H$. Theorem~\ref{T6}$'$ shows that for a given $R$-matrix there is a 1-1 correspondence between
$*$-canonical elements for a $*$-QTQHA $H$ and quasi-triangular compatible twists on $H$.

It is worth noting that
\begin{Lemma} Let $(H,\O,\R)$ be a $*$-QTQHA. Then $H$ is also a $*$-QTQHA with $*$-canonic\-al element
$\O^\dag$.
\end{Lemma}
\begin{proof}
       Taking the conjugate inverse of equation (\ref{E39}) gives
       $$\R={(\O^T)^{-1}}^\dag(\R^\dag)^{-1}\O^\dag\quad[\textrm{resp.} {(\O^T)^{-1}}^\dag (\R^T)^\dag \O^\dag]$$
       or equivalently
       $$(\R^\dag)^{-1}=(\O^\dag)^T \R (\O^\dag)^{-1}\quad[\textrm{resp.} (\O^\dag)^T (\R^T)^{-1} (\O^\dag)^{-1}]$$
       which, together with proposition~\ref{P10}, is sufficient to prove the result.
\end{proof}
\begin{Corollary*} $\O^{-1}\O^\dag$ must determine a quasi-triangular compatible twist on $H$.
\end{Corollary*}

This last result puts a strong restriction on $\O$ in order for it to give rise
to a $*$-canonical element for a $*$-QTQHA.

\newpage

\end{document}